\documentclass[10pt,reqno]{amsart} 			%% texte 
\usepackage{amsmath, amscd, amsthm, amsfonts, amssymb}
\usepackage{mathrsfs, bbm, cancel,  stmaryrd, xfrac}
\usepackage{enumerate, enumitem, caption}
\usepackage{slantsc, lmodern, comment}
\usepackage[all]{xy}
\usepackage{graphicx}
\usepackage{rotating, lscape, longtable, multirow, array, cite}
\usepackage{titletoc}
\usepackage{dsfont}

\usepackage{color}
\usepackage{hyperref}
\definecolor{darkblue}{rgb}{0,0,1}
\definecolor{darkgreen}{rgb}{0,1,0}
\hypersetup{
pdftitle={11D SUSY CW-space},
pdfauthor={Frank Klinker},
pdfkeywords={supersymmetry, CW-space, superalgebra},
pdfpagemode=UseNone,
linktocpage=true,
breaklinks=true,  
colorlinks=true,
linkcolor=darkblue,
urlcolor=darkblue, 
citecolor=darkblue%red
}

\newtheorem{prop}{Proposition}[section]
\newtheorem{thm}[prop]{Theorem}

\theoremstyle{definition}

\newtheorem{defn}[prop]{Definition}

\newtheorem{rem}[prop]{Remark}
\newtheorem{exmp}[prop]{Example}

\theoremstyle{remark}

\newcommand{\cl}{{\rm C\ell}}

\newcommand{\q}{{\rm q}}

\newcommand{\ZZ}{\mathbbm Z}
\newcommand{\RR}{\mathbbm R}

\newcommand{\CC}{\mathbbm C}

\newcommand{\ei}{\mathbbm{1}}

\newcommand{\eucl}{\langle\cdot,\cdot\rangle}

\allowdisplaybreaks

\renewcommand{\tabcolsep}{2pt}
\renewcommand{\arraystretch}{1}

\begin{document}

\title{A family of non-restricted $D=11$ geometric supersymmetries}
\author{Frank Klinker}
\thanks{\hspace*{-1em}2000 {\em Mathematics Subject Classification}: 
17B66, 17B81, 53B30, 53C27, 83E50}
\thanks{\hspace*{-1em}{\em Key words and phrases}: 
Geometric supersymmetry, geometric superalgebra, spinor connection, Cahen-Wallach space, supergravity background}
\thanks{\ \\[-2ex]\hspace*{0.5em}{\em Published in}:
\href{http://dx.doi.org/10.1016/j.geomphys.2014.10.008}{J.~Geom.~Phys.~{\bf 86} (2014) 534-553}}
%%\thanks{online: 2012-xx-xx.}
\thanks{\ \\[-2ex]\hspace*{0.5em}{\bf Note}: 
Compared with the published version we clarified the content of Section \ref{sec:7} and corrected one of the interpretations.}
\address{
Faculty of Mathematics, TU Dortmund University, 44221 Dortmund, Germany
} 
\email{
\url{frank.klinker@math.tu-dortmund.de}
}
\begin{abstract} 
We construct a two parameter family of eleven-dimensional indecomposable Cahen-Wallach spaces with irreducible, non-flat, non-restricted geometric supersymmetry of fraction $\nu=\sfrac{3}{4}$. 
Its compactified moduli space can be parametrized by a compact interval with two points corresponding to two non-isometric, decomposable spaces. 
These singular spaces are associated to a restricted $N=4$ geometric supersymmetry with $\nu=\sfrac{1}{2}$ in dimension six and a non-restricted $N=2$ geometric supersymmetry with $\nu=\sfrac{3}{4}$ in dimension nine.
\end{abstract}

\maketitle
\setcounter{tocdepth}{1}
\tableofcontents 

\section{Introduction}

In this text we describe in detail the geometric supersymmetry of a family of eleven dimensional manifolds. Geometric supersymmetry is by definition an extension of the Lie algebra of Killing vector fields to a super Lie algebra by geometric data. This roughly means that the odd part of the superalgebra is given by a linear subspace of the sections in a bundle over the manifold compatible with the Killing vector fields -- we will be more detailed in the beginning of Section \ref{sec:geoalg}. 
Although the manifolds we consider in this text are homogeneous spaces and, therefore, it would be sufficient to discuss the structure at one point, we will give all results in terms of local coordinates.
On the one hand, we do this to emphasize the geometric nature and, on the other hand, because the local description may give an idea for similar constructions in the non homogeneous situation. However, the calculations are similar in both concepts.

First we will provide the general setup of the manifolds that we consider, the Cahen-Wallach spaces. We describe the local structure of the metric and determine the Killing vector fields that yield the even part of our structure, see Sections \ref{sec:gen} and \ref{sec:metric+Killing}.
Then we turn to the odd part: It will be spanned by sections in a spinor bundle that are parallel with respect to a given connection. Again, we will give the local description and show how this depends on the elements at one point, see Section \ref{sec:oddPart}.
These preparations lead to Section \ref{sec:geoalg} where we first give a short introduction to geometric superalgebras and geometric supersymmetries. Then we describe in detail how the ingredients provided so far define a geometric superalgebra. In particular, we prove several compatibility conditions. 
In Section \ref{sec:susy} we discuss whether there are situations in which the geometric superalgebra yields geometric supersymmetry. We formulate the obstruction and provide a full list. 
In the final Section we discuss the moduli space of geometric superalgebras  and geometric supersymmetries. Furthermore, we associate the singularities to extended geometric supersymmetries in dimensions six and nine.

We would like to mention that the spaces we discuss in this text are canonical candidates for supergravity backgrounds. For more details on the supergravity point of view that go beyond the explanations in our text we cordially refer the reader to the literature in the references, e.g. \cite{Pope:2002,FigPapado1,Gauntlett:2002,MeFig04,Hustler,Meessen2002}.

\section{The general setup}\label{sec:gen}

The classification of solvable Lorentzian symmetric spaces by the construction presented below goes back to \cite{CW70}. 
Let $(V,\eucl)$ be an $n$-dimensional euclidean vector space and $B$ be a symmetric endomorphism of $V$. We denote the symmetric bilinear form that is defined by $B$ and $\eucl$ by the same symbol $B$ and we write $*:V\to V^*$, $v\mapsto v^*$ with $v^*(w)=\langle v,w\rangle$ for the canonical identification of $V$ and its dual. We define $W:=\RR^{1,1}\oplus V$ and denote by $\mathring{g}$ the extension of $\eucl$ to a block diagonal Lorentzian metric on $W$. We consider a null basis  $\{e_+,e_-\}$ of $\RR^{1,1}$ with respect to $\mathring{g}|_{\RR^{1,1}}$. The following skew symmetric multiplication on $\mathfrak{g}:=V^*\oplus W$ yields a Lie algebra structure on $\mathfrak{g}$:
\begin{align}
\left[e_-,w\right] &= w^*\,, \label{c1} \\
\left[v^*, e_-\right] &= Bv\,, \label{c2}\\
\left[v^*,w\right] &= -v^*(Bw)\cdot e_+ =  -\langle Bv,w \rangle\cdot e_+\,, \label{c3}
\end{align}
for all $w\in V$ and $v^*\in V^*$. The bilinear form $\mathring{g}$ is extended to a bi-invariant metric on $\mathfrak{g}$ by $\mathring{g}(v^*,w^*):=\langle Bv,w\rangle$ if $B$ is non degenerate, see below.

Within $\mathfrak{g}$ the factor $V^*$ acts on $W$, 
the bracket of $W$ with itself obeys $[W,W]= V^*$, 
and $\mathring{g}$ is $V^*$-invariant. 
From (\ref{c1})-(\ref{c3}) we see that the embedding 
\begin{equation}
V^*\longrightarrow \RR_+\otimes V \hookrightarrow \mathfrak{so}(W)
=\mathfrak{so}(V)\,\oplus\, (\RR_+\otimes V) \,\oplus\,(\RR_-\otimes V)\,\oplus\,( \RR_+\otimes\RR_-)
\end{equation}
is given by $v^*\mapsto  Bv\wedge e_+$ where $x\wedge y (z):=\langle y,z\rangle x-\langle x,z\rangle y$.

The above data yield a ($D=n+2$)-dimensional symmetric space $M_B$ with Lorentz metric determined by $\langle\cdot,\cdot\rangle$ and $B$. 
The resulting Lorentzian space $M_B$ is indecomposable if and only if the symmetric map $B$ is non-degenerate. This can best be seen from (\ref{c1})-(\ref{c3}) if we recall that $M_B$ is decomposable if there exists a $V^*$-invariant subspace $\tilde W\subset W$ such that $\mathring{g}|_{\tilde W\times \tilde W}$ is non-degenerate, see \cite{Baum12, Wu67}. 
If $B$ admits zero eigenvalues $\mathring{g}$ is degenerate but it remains a metric if $V^*$ is truncated. In this case the resulting manifold decomposes into a product of a lower dimensional Cahen-Wallach space and an euclidean space.\footnote{If we refer to $\mathring{g}$ as metric we will always assume this truncation.} This can also be deduced from the coordinate form of the metric, see (\ref{11}).

\section{The metric and the Killing vector fields}\label{sec:metric+Killing}

\subsection{The metric}
For the coordinate description of $M_B$ we may use the exponential map and write for $x=x^+e_++x^-e_-+\vec{x}\in W$ with $\vec x=\sum_ix^ie_i$ 
\[
\mu(x):=\exp\big(x^+e_+\big)\exp\big(x^-e_-\big)\exp\big(\vec x\big)\,.
\]
This obeys
\begin{align}
\partial_+\mu &= \exp(x^+e_+)e_+ \exp(x^-e_-)\exp(\vec x)
			  =  \mu(x) e_+ 
\\			  
\partial_i\mu &= \exp(x^+e_+)\exp(x^-e_-)\exp(\vec x) e_i
			  = \mu(x) e_i
\\				  
\partial_-\mu &= \exp(x^+e_+)\exp(x^-e_-) e_- \exp(\vec x)	
			  = \mu(x)\exp(-\vec x ) e_-\exp(\vec x)\nonumber\\
			  & = \mu(x)\textstyle\big(e_-+\sum_ix^ie_i^*- \tfrac{1}{2}\sum_{ij}B_{ij}x^ix^j e_+\big)
\end{align}
where we use $\exp(\vec x)=\prod_i\exp(x^ie_i)$ and 
\begin{equation}
\begin{aligned}
e_j^*\exp(x^ie_i)&=\exp(x^ie_i)\big(e_j^*-B_{ij}x^ie_+\big)\\
e_- \exp(x^ie_i)&= \exp(x^ie_i)\big( e_- + x^ie_i^* -\tfrac{1}{2}B_{ii}(x^i)^2 e_+\big)
\end{aligned}
\end{equation}
with the symmetric matrix $(B_{ij})$ defined by $B(e_i)=\sum_jB_{ji}e_j$.  

From this we read the two components of the Maurer-Cartan form 
$\mu^{-1}d\mu = \omega+\theta \in \Omega^1(M_B)\otimes\mathfrak{g}$ 
with $\omega\in\Omega^1(M_B)\otimes V^*$ and 
$\theta\in\Omega^1(M_B)\,\otimes\, W$:
\begin{align}
\omega & = \sum_i x^i dx^-\otimes e^*_i\,, \\ 
\theta & = dx^-\otimes e_- +\sum_i dx^i\otimes e_i
		+\big(dx^+-\tfrac{1}{2}\sum_{ij} B_{ij}x^i x^j dx^-\big)\otimes e_+\,.
\end{align}  

With $g_B= \mathring{g}(\theta,\theta)$ we get the following local form of the metric on $M_B$:
\begin{equation}\label{11}
g_B  = 2dx^+dx^- - \sum_{ij} B_{ij}x^i x^j(dx^-)^2 + \sum_i (dx^i)^2\,.
\end{equation} 

In particular, the Levi-Civita connection of $g_B$ is determined by the Christoffel symbols 
$\Gamma_{i-;-}=-\Gamma_{--;i}= -\sum_jB_{ij}x^j$. 
If we move from the coordinates to the adapted ON frame $\{\partial_+,\partial_-+\frac{1}{2}\sum_{ij}B_{ij}x^ix^j\partial_+,\partial_i\}$  there is only one surviving component of the connection form, namely
\begin{equation}
\omega_{i-}=-\omega_{-i}=\sum_{j}B_{ij}x^jdx^-\,.
\end{equation}

The bi-invariant metric $\mathring{g}$ on $\mathfrak{g}$ makes the decomposition $\mathfrak{g}=V^*\oplus W$ an orthogonal splitting and the isometry algebra of $M_B$ is given by
\begin{equation}
\mathfrak{isom}(M_B)= \mathfrak{so}_B(V)\oplus V^*\oplus W
\end{equation}  
with 
\begin{equation}\label{condition:soB}
\begin{aligned}
\mathfrak{so}_B(V)&=\big\{A\in\mathfrak{so}(V)\,|\,[A,B]=0\big\}\\
&=\big\{A\in\mathfrak{so}(V)\,|\, [A,v^*]=(Av)^*\text{ for all }v^*\in V^*\big\}\,.
\end{aligned}
\end{equation}
\begin{rem}\label{rem:conformal}
To fix some notation, we like to mention that two Lorentzian spaces defined by symmetric maps $B_1$ and $B_2$ are isometric if and only if $B_1$ and $B_2$ are conformally equivalent, i.e.\ there exists a real scalar $c>0$ and an orthogonal transformation $X$ such that $B_2=c X^tB_1 X$. 

Therefore, we assume $B$  to be diagonal such that the space $M_B$ is defined by a sequence of real numbers $\lambda_1^2,\ldots,\lambda^2_n$. 
We may also sort this sequence in the way $\lambda_1^2\leq \ldots\leq \lambda_{n}^2$ and all non-vanishing if $M_B$ is indecomposable. We write the eigenvalues as squares of real or imaginary numbers because their square roots will play an important role in our calculations.
\end{rem}

\begin{exmp}\label{ex:maximal}
Consider $D=11$, i.e.\ $n=9$, and $B=-4\beta^2\begin{pmatrix}4\mathbbm{1}_3 & \\ & \mathbbm{1}_6\end{pmatrix}$. Then $M_B$ is indecomposable and the metric is given by 
\[
g_B=2dx^+dx^- +4\beta^2 \Big(4\sum\limits_{i=1}^3(x^i)^2+\sum\limits_{i=4}^9 (x^i)^2\Big)(dx^-)^2 + \sum\limits_{i=1}^9 (dx^i)^2\,.
\]
%This particular example has been considered in \cite{CheKo84} and \cite{FigPapado1} as maximally symmetric background of eleven dimensional supergravity.
\end{exmp}

\subsection{The Killing vector fields}\label{sec:Killing}

A local basis of the isometry algebra of $M_B$ is provided by the Killing vector fields, i.e. by those vector fields $X$ that obey $L_Xg_B=0$. We will denote the Killing vector fields associated to the ON frame of $\mathfrak{g}$ by $K_{(+)},K_{(-)},K_{(i)},K_{(i^*)}$ and those associated to the standard basis of $\mathfrak{so}_B(V)$ by $K_{(ij)}$. From here we consider $B={\rm diag}(\lambda_1^2,\ldots,\lambda^2_n)$.

Because the metric coefficients only depend on the $x^i$ we immediately see that $\partial_+$ and $\partial_-$ are Killing vector fields and we write $K_{(+)}=-\partial_+$ and $K_{(-)}=-\partial_-$. 
The ansatz 
\begin{equation}
\begin{aligned}
K_{(i)}&=\alpha_i(x^-)\partial_i+\beta_i(x^-)x^i\partial_+\,,\\
K_{(i^*)}&=\alpha^*_i(x^-)\partial_i+\beta^*_i(x^-)x^i\partial_+\,,
\end{aligned}
\end{equation}
inserted into $\mathcal{L}_Kg_B=0$ yields 
\[
\frac{\partial \beta^{(*)}_i}{\partial x^-} = \lambda^2_i\alpha^{(*)}_i\,, \  
\frac{\partial\alpha^{(*)}_i}{\partial x^-} = -\beta^{(*)}_i \,,
\]
or 
\[
\frac{\partial^2 \beta^{(*)}_i}{\partial (x^-)^2} = -\lambda^2_i\beta^{(*)}_i\,, \  
\frac{\partial^2\alpha^{(*)}_i}{\partial (x^-)^2} = -\lambda^2_i\alpha^{(*)}_i \,.
\]
This motivates the further specialization to 
\[
\alpha_i= a_i\cos(\lambda_ix^-)\,,\  \beta_i=b_i\sin(\lambda_ix^-)\,,\ 
\alpha_i^*= a^*_i\sin(\lambda_ix^-)\,,\  \beta^*_i=b^*_i\cos(\lambda_ix^-)\,,
\]
and the coefficients are related by $\lambda_ia_i=b_i$, $-\lambda_ia_i^*=b^*_i$. 
By claiming the commutation relations (\ref{c1})-(\ref{c3}) we fix the remaining free parameters:\footnote{We have to take into account that the vector fields obey the commutation relations only up to sign. This is due to the difference between right and left invariance when we turn from the group structure to the structure on the coset space, see \cite{KobNom} for more details on this fact.} (\ref{c1}) and (\ref{c2}) yield $a_i^*=-\lambda_ia_i$ and (\ref{c3}) yields $a_i^2=1$. The Killing vector fields that are adapted to $e_+,e_-$ and the orthonormal eigenbasis $\{e_i\}$ of $B$ as well as their duals $\{e_i^*\}$ are $K_{(+)}=-\partial_+$, $K_{(-)}=-\partial_-$, and 
\begin{align}
K_{(i)}&=	\cos(\lambda_ix^-)\partial_i 
		+ \lambda_i\sin(\lambda_ix^-)x^i\partial_+\,, \label{Ki}\\
K_{(i^*)}&=	-\lambda_i\sin(\lambda_ix^-)\partial_i
		+\lambda_i^2 \cos(\lambda_ix^-)x^i\partial_+\,.\label{Ki*}
\end{align}
The additional Killing vector fields that come from $\mathfrak{so}_B(V)$ are given by the usual $\mathfrak{so}(V)$ generators subject to condition (\ref{condition:soB}), i.e.
\begin{equation}\label{Kij}
K_{(ij)}=x^j\partial_i-x^i\partial_j
\end{equation}
with $i,j\in I_\alpha$ for some $\alpha$ where $\{1,\ldots,n\}=\bigcup_{\alpha=1}^{\hat n} I_\alpha$ with $I_\alpha=\{r_{\alpha-1}+1,\ldots,r_\alpha\}$ is the decomposition coming from $\lambda_1^2=\ldots=\lambda^2_{r_1}<\lambda^2_{r_1+1}=\ldots=\lambda^2_{r_2}<\ldots<\lambda^2_{r_{\hat n-1}+1}=\ldots=\lambda^2_{n}$ and $r_0=1,r_{\hat n}=n$.

\section{Connections and parallel spinors}\label{sec:oddPart}

We will consider a special class of connections on a spinor bundle $S(M_B)$ of $M_B$, namely connections that are compatible with the homogeneous structure of $M_B$. 
If $S(M_B)$ is associated to the irreducible Clifford module such connections are described by $V^*$-equivariant linear maps $\mathfrak{g}=V^*\oplus W\to \cl(W)$ with the property that $\rho:\mathfrak{so}(W)\supset V^* \to\cl(W)$ coincides with the spin representation, i.e.\ $\rho(v^*)=\Gamma(v^*)$. Here $\cl(W)$ denotes the Clifford algebra of $W$. In \cite{klinker:CW} such connections are discussed in detail, and we will state the result in Propositions \ref{prop:maps} and \ref{prop:flat}. In Appendix \ref{sec:conv} we recall our conventions for the Clifford algebra. Moreover we will collect some facts that will be frequently used in our text, in particular, the calculations from Section \ref{sec:geoalg} will become more transparent.

\subsection{Preliminaries}\label{sec:prel}
Mainly to fix our notation, we recall some facts on the Clifford algebra in the case of Cahen-Wallach spaces. We consider $\cl(\mathbbm{R}^{1,1})=\mathfrak{gl}_2\CC$ with generators 
$\gamma_+=\gamma(e_+)=\frac{1}{\sqrt{2}}(i\sigma_2+\sigma_1)
		=\sqrt{2}{\begin{pmatrix}0&1\\0&0\end{pmatrix}}$, $
 \gamma_-=\gamma(e_-)=\frac{1}{\sqrt{2}}(i\sigma_2-\sigma_1)
 		=\sqrt{2}{\begin{pmatrix}0&0\\-1&0\end{pmatrix}}$ 
and we denote the two-dimensional volume element by 
$\sigma:=\frac{1}{2}\left[\gamma_+,\gamma_-\right]=-\sigma_3={\begin{pmatrix}-1&0\\0&1\end{pmatrix}}$. 
If we denote the generators of $\cl(V)$ by $\{\gamma_i\}_{1\leq i\leq n}$ those of 
$\cl(W)=\mathfrak{gl}_2\CC \,\hat\otimes \,\cl(V)$  are given by 
\begin{equation}\label{basis}
\{\Gamma_\mu\}_{\mu\in\{+,-,i\}}=\{\gamma_+\otimes\ei, \gamma_{-}\otimes\mathbbm{1},\sigma\otimes\gamma_i\}\,.
\end{equation} 
In particular, 
\begin{align}
\mathfrak{gl}_2\CC \ni r &\mapsto r\,\hat\otimes\,\ei = r \otimes\mathbbm{1} \in \cl(W)\,,\label{incl-gl}\\
\cl(V)\ni a &\mapsto \mathbbm{1}\,\hat\otimes\,a=\mathbbm{1}\otimes a^0+ \sigma\otimes a^1\in \cl(W)\,,\label{incl-cl}
\end{align}
where $a=a^0+a^1\in \cl(V)$ is the decomposition into its even and odd part. In this regard, we consider the map $\bar{\ }:\cl(V)\to\cl(V)$ with $\overline{a^0+a^1}=a^0-a^1$. By $\hat\otimes $ we denote the $\ZZ_2$-graded tensor product that describes the relation between the Clifford algebra of the sum of vector spaces and the Clifford algebras of the respective summands, see \cite[Prop. 1.5]{LawMich}. In terms of the usual tensor product the graded isomorphism is achieved by introducing the volume form in (\ref{basis}).
 
Consider the irreducible Clifford modules $S_2$ and $S(V)$ of $\cl(\mathbbm{R}^{1,1})$ and $\cl(V)$.
The first one decomposes into a sum of two one dimensional half spinor spaces $S_2^\pm={\rm ker}(\gamma_\mp)$ given by the  $\pm1$-eigenspaces of $\sigma$. If we denote the two projections on the two eigenspaces by $\sigma_\pm=\frac{1}{2}(\ei\pm \sigma)=-\frac{1}{2}\gamma_{\mp}\gamma_{\pm}$ then (\ref{incl-cl}) is rewritten as 
$\ei\,\hat\otimes\, a=\sigma_-\,\hat\otimes\,a+\sigma_+\hat \otimes\,a=\sigma_-\otimes \bar a+\sigma_+\otimes a$.
In our choice of $\gamma$-matrices the eigendirections are given by $\vec{e}_1=(1,0)^t$ and $\vec{e}_2=(0,1)^t$ such that a spinor in $S(W)= S_2\,\hat\otimes\,S(V)= S_2^-\otimes S(V)\oplus S^+_2\otimes S(V)=:S^-(W)\oplus S^+(W)$ can be written as 
$\vec\eta={\renewcommand{\arraystretch}{1.3}\begin{pmatrix}\eta_1\\\eta_2\end{pmatrix}}= \vec{e}_1\otimes\eta_1+ \vec{e_2}\otimes \eta_2$. 
The action of $\cl(W)$ on $S(W)$ is given by 
\begin{equation}\label{action:wichtig}
\begin{aligned}
{(r \,\hat\otimes\,a)\begin{pmatrix}\eta_1\\\eta_2 \end{pmatrix} 
=\begin{pmatrix}r_{11}\bar a & r_{12}a\\r_{21} \bar a & r_{22}a\end{pmatrix}
	\begin{pmatrix}\eta_1\\\eta_2 \end{pmatrix} 
= r\begin{pmatrix}\bar a\eta_1\\ a\eta_2
\end{pmatrix}}
\end{aligned}
\end{equation}
for $r\in\mathfrak{gl}_2\CC$ and $a\in\cl(V)$. 
In particular, the image of $v^*\in V^*$ considered as an element of $\mathfrak{so}(W)\subset \cl(W)$ under the spin representation is given by 
\begin{equation}\label{vstar}\begin{aligned}
v^*= Bv\wedge e_+ \mapsto \ &
	\frac{1}{4}\big(( \gamma_+\otimes\ei)( \sigma\otimes Bv)-(\sigma\otimes Bv)( \gamma_+ \otimes\ei)\big) \\
&=	\frac{1}{4} (\gamma_+\sigma-\sigma\gamma_+)\otimes  Bv 
= \frac{1}{2}\gamma_+\otimes Bv
= \frac{1}{\sqrt{2}}{\begin{pmatrix}0&Bv\\0&0\end{pmatrix}} \,,
\end{aligned}\end{equation}
see Appendix \ref{sec:conv}.
\subsection{The algebraic description}

In terms of the notation introduced above the relevant spinor connections of $M_B$ are singled out as follows.

\begin{prop}\label{prop:maps}
The $V^*$ equivariant linear maps that define homogeneous connections on the spinor bundle are 
\[
\begin{aligned}
\rho(v^*) &= \frac{1}{2}\gamma_+\hat\otimes Bv 
		= \frac{1}{\sqrt{2}}{\begin{pmatrix}0&Bv\\0&0\end{pmatrix}}\,,\\
\rho(e_+) &= \frac{1}{2}\gamma_+\hat\otimes a 
		= {\begin{pmatrix}0& \sqrt{2}\,a\\ 0 &0\end{pmatrix}}\,,\\
\rho(e_-) &= \sigma_-\hat\otimes c +\sigma_+\hat\otimes d+\gamma_-\hat\otimes b +\gamma_+\hat\otimes \epsilon
		= {\begin{pmatrix} \bar c & \sqrt{2}\, \epsilon \\ \sqrt{2}\, b & d\end{pmatrix}}\,,\\
\rho(w)   &= -\sigma_-\hat\otimes wb -\sigma_+\hat\otimes bw - \frac{1}{2} \gamma_+\hat\otimes s_{\bar c,d}(w)
		={\begin{pmatrix}w  b & -\frac{1}{\sqrt{2}}s_{\bar c,d}(w)\\ 0 & -bw\end{pmatrix}}\,,
\end{aligned}
\]
with $a,b,c,d,\epsilon\in\cl(V)$ and 
\[
s_{\bar c,d}:\cl(V)\to \cl(V),\quad s_{\bar c,d}(x)=\bar cx-xd\,.
\]
The two parameters $a,b$ are fixed to be pseudo-scalars $a=\alpha+\beta\gamma^*$ and $b= -\alpha+\beta\gamma^*$ if ${\rm dim}(V)$ is even, and scalars $a=-b= \alpha$ if ${\rm dim}(V)$ is odd.\footnote{In particular, $b$ is an even Clifford element such that an additional $\bar{\ }$ is not needed when we turn to the matrix description.}
\end{prop}

\begin{rem}\label{rem:soBV}
We consider $\mathfrak{so}_B(V)$ acting in the usual way on $W$. Then it is compatible with the equivariant map $\rho$ if it is extended by $\rho(A):=\Gamma(A)$ for all $A\in\mathfrak{so}_B(V)$ where $\Gamma$ is the spin representation.

As we will see in Proposition \ref{prop:even-odd} and Remark \ref{rem:liederivative} the map $\rho$ and its extension to $\mathfrak{so}_B(V)$ will play a crucial role when we consider the action of the Killing vector fields on the spinors. In particular, we know that there exist additional Killing vector fields that are not directly connected to the Lie algebra $\mathfrak{g}$ but to $\mathfrak{so}_B(V)$, namely (\ref{Kij}).
\end{rem}

The curvature of a connection given by the equivariant map $\rho$ is determined by its values on $W$ and given by
\[
\mathcal{R}^\rho(X,Y)=[\rho(X),\rho(Y)] - \rho([X,Y]_W) -\Gamma([X,Y]_{V^*})\,.
\]
Therefore, an equivariant map $\rho$ from Proposition \ref{prop:maps} yields a flat connection if and only if $\rho$ is a representation. 

If we assume scalar parameters $a=-b=\alpha$, for example, the surviving components of the curvature are
\begin{align}
R^\rho(e_-,e_+)
	&=\begin{pmatrix} 
		-\alpha^2 & \frac{\alpha}{\sqrt{2}}(\bar c-d)\\ 
		0 & \alpha^2
	  \end{pmatrix}\,,\\
R^\rho(e_i,e_j) 
	&=\begin{pmatrix}
		2\alpha^2\gamma_{ij} & -\sqrt{2}\alpha\{s_{\bar c,d}(e_{[i}),\gamma_{j]}\}\\
		0 & 2\alpha^2\gamma_{ij}
	  \end{pmatrix}\,,\label{curv1}\\
R^\rho(e_-,e_i) 
	&=-\frac{1}{\sqrt{2}}\begin{pmatrix}	
			0& q_{\bar c,d}(e_i)+ B(e_i) \\	
			0& 0
	   \end{pmatrix}\,,\label{curv-MB}
\end{align}
with 
\begin{equation}
q_{\bar c,d}:\cl(V)\to\cl(V),\quad q_{\bar c,d}(x)=s_{\bar c,d}^2(x)=\bar c^2x+xd^2-2\bar cxd\,.
\end{equation}

The flat connections are singled out from Proposition \ref{prop:maps} by the following Proposition \ref{prop:flat} along with Remark \ref{rem:e+}.
\begin{prop}\label{prop:flat}
An equivariant map $\rho$ with $\rho(e_+)=0$ defines a flat connection if and only if
\begin{equation}\label{rep-alpha0}
\begin{gathered}
\rho(v^*) = \frac{1}{\sqrt{2}}\begin{pmatrix}0&Bv\\0&0\end{pmatrix}\,,\quad
\rho(w)   = -\frac{1}{\sqrt{2}}\begin{pmatrix}0& s_{\bar c,d}(w)\\ 0 &0\end{pmatrix}\,,\\
\rho(e_-) = \begin{pmatrix}\bar c & \sqrt{2}\, \epsilon \\ 0 & d\end{pmatrix}\,,
\end{gathered}
\end{equation}
with $(\bar c,d)$ subject to
\begin{equation}\label{qcp}
q_{\bar c,d}(v)=-B(v)
\end{equation}
for all $v\in V$. In particular, $\epsilon\in \cl(V)$ remains a free parameter of the flat connection.
\end{prop}

\begin{exmp}[Example \ref{ex:maximal} continued]\label{ex:maximal-cont}\ For the metric associated to the symmetric map $B=-4\beta^2\begin{pmatrix}
4\mathbbm{1}_3 &\\&\mathbbm{1}_6\end{pmatrix}$ in eleven dimensions, the pair $(\bar c,d)$ with
\[
\bar c=-3\beta\Gamma_{123}\,,\quad d=\beta\Gamma_{123}\,,
\]
obeys (\ref{qcp}). This means, the spinor connection that is defined by these data is flat. In fact this pair along with $\epsilon=0$ has been considered in \cite{CheKo84,FigPapado1} as a  connection that provides a maximal amount of parallel spinors.
\end{exmp}

In \cite{klinker:CW} we discuss in detail a large class of pairs $(\bar c,d)$ that solve condition (\ref{qcp}), the so called quadratic Clifford pairs. Furthermore we give an additional condition that makes the list of solutions we present complete. This condition arises naturally in the discussion of supersymmetry.
 
Although we will not need it later on, but for the sake of completeness, we will state the result analog to Proposition \ref{prop:flat} for $\rho(e_+)\neq 0$. 
\begin{rem}\label{rem:e+}
The equivariant map $\rho$ with $\rho(e_+)\neq 0$ defines a flat connection only if $n$ is even and  $B=-2\lambda^2\ei$. Moreover, $a=\Pi^\pm$ is required. We consider the upper sign and the map $\rho$ is given by
\begin{equation}\label{rep-alphanot0}
\begin{aligned}
\rho(v^*) 	=\ & -\sqrt{2} \lambda^2 \begin{pmatrix} 0& 	v \\ 0& 0 \end{pmatrix}\,,\qquad
\rho(e_+) 	=\  \sqrt{2}\alpha\begin{pmatrix} 0& \Pi^+ 	 \\ 0& 0 \end{pmatrix} \,,
\\
\rho(e_-) 	=\ & \sqrt{2}\begin{pmatrix}
					\rho_0 - \sqrt{\alpha\beta+\lambda^2} - \bar c^+_- &\beta\Pi^- +(\epsilon^+_-+\epsilon^-_++\epsilon^+_+)  \\
					  -\alpha\Pi^-  & \rho_0+\sqrt{\alpha\beta+\lambda^2} + d^-_+
				 \end{pmatrix}\,, \\
\rho(v)  	=\ & \begin{pmatrix}
			    -\alpha \Pi^+ v  & \sqrt{\alpha\beta+\lambda^2} v + s_{\bar c_-^+,d_+^-}(v) \\
			   0  &  \alpha\Pi^- v
			\end{pmatrix}\,.
\end{aligned}
\end{equation}
The free parameters are the scalars $\alpha,\beta,\rho_0$ and the Clifford element $\epsilon^+_+$. The further contributions are related by 
$\sqrt{\alpha\beta+\lambda^2}\,s_{\bar c^+_-,d^-_+}(v)=\sqrt{2}\alpha\, s_{\epsilon^+_-,-\epsilon^-_+}(v)$ for all $v\in V$.
\end{rem}

\subsection{The local description}\label{subsec:local}

The discussion so far took place in one particular point of the manifold and, therefore, was a discussion on the level of Lie algebras and representations. 
Below we will use the coordinates introduced above to get a local description on the manifold $M_B$.

The choice of coordinates on $M_B$ yields a splitting of the spinor bundle $S:=S(M_B)$ of $M_B$ as $S^{-}\oplus S^{+}$ with the first (resp.\ second) summand
being the $-1$-eigenspaces (resp.\ $+1$-eigenspace) of $\sigma:=\frac{1}{2}[\Gamma_+,\Gamma_-]$. 
The projections on the two subbundles are given by $\sigma_{\pm}=-\frac{1}{2}\Gamma_{\mp}\Gamma_{\pm}$.    
We will use here the notation $\vec\xi=\xi_1+\xi_2$ for the sections in $S$, too. 
Due to $\Gamma_+^2=0$ we have $\Gamma_+:S^-\to S^+$ and $S^-={\rm ker}(\Gamma_+)$. Furthermore, both subbundles are preserved by the action of $\Gamma_i$. 

The Levi-Civita connection on $M_B$ induces a connection on the spinor bundle $S$ via the spin representation. 
It is given by 
$\nabla\vec\xi=d\vec\xi-\frac{1}{4}\sum_{\mu\nu}\omega^{\mu\nu}\Gamma_{\mu\nu}\vec\xi$ 
which in our situation is
\begin{equation}
\begin{aligned}
\nabla_+\vec\xi
	&= \partial_+\vec\xi\,, \\
\nabla_-\vec\xi
	&= \partial_-\vec\xi -\frac{1}{2}\omega_-^{+i}\Gamma_{+i} \vec\xi
	 = \partial_-\vec\xi +\frac{1}{2}\sum_ix^i\Gamma_+ B(e_i)\xi_2\,,\\
\nabla_i\vec\xi
	&= \partial_i\vec\xi \,.
\end{aligned}
\end{equation}

We will consider a general connection on $S(M_B)$ defined by the equivariant map $\rho$ as given in Proposition \ref{prop:maps}. 
For the discussion of spinor connections we restrict ourselves to the situation of Proposition \ref{prop:maps} with scalars $a=-b=\alpha$ because we will discuss odd dimensional manifolds later on. 
We know about the dimension of the space of parallel sections, $\mathcal{K}_1\subset\cancel S$: It coincides with the dimension of the kernel of the curvature $R^\rho$. For a connection according to Proposition \ref{prop:flat} we get $\dim \mathcal{K}_1=\dim S(W)={\rm rank}\,S$. More general, we see from (\ref{curv1}) and (\ref{curv-MB}) that $\alpha=0$ is required if we assume the kernel of $R^\rho$ to be non-trivial. Moreover, in this case the kernel is given by ${\rm ker}(\Gamma_+)= S^-$, generically. 

We will consider this situation with $\epsilon=0$ such that the connection is entirely determined by the pair $(\bar c,d)$. In our local coordinates the connection is given by 
$D_\mu=\nabla_\mu+\rho(e_\mu)$ for $\mu\in\{+,-,i\}$ and the parallel spinors satisfy $D\vec\xi=0$.
After applying the projecton operators $-\frac{1}{2}\Gamma_\pm\Gamma_\mp$ this is 
\begin{align}
0&= \partial_+\xi_\alpha\,, \label{f1} \\
0&=\partial_i\xi_1 -\frac{1}{2}\Gamma_+s_{\bar c,d}(e_i)\xi_2\,,\label{f2}\\
0&=\partial_i\xi_2\,, \label{f3}\\
0&=\partial_-\xi_1 + \frac{1}{2}\Gamma_+\sum_jx^jB(e_j)\xi_2 + \bar c\xi_1\,, \label{f4}\\
0&=\partial_-\xi_2 + d\xi_2\,.\label{f5}
\end{align}
From (\ref{f1}) we see that $\vec\xi$ is independent of $x^+$ and from (\ref{f3}) that $\xi_2$ is independent form $x^i$. Therefore, (\ref{f5}) or $\partial_-\xi_2 =- d\xi_2$ yields
\begin{equation}
\xi_2=\xi_2(x^-)=\exp(-x^- d)\xi_2^0\label{C}
\end{equation}
for a constant spinor $\xi^0_2$.

Moreover, (\ref{f2}) yields $\partial_i\partial_j\xi_1=0$ such that
\begin{equation}\label{D}
\xi_1=\xi_1(x^-,x^j)
=\xi'_1(x^-)+\frac{1}{2}\sum_i x^i\Gamma_+s_{\bar c,d}(e_i)\xi_2(x^-)
\end{equation}
for a spinor $\xi'_1$ depending only on $x^-$.

Inserting both in (\ref{f4}) yields
\[
\begin{aligned}
0&=
\partial_-\xi'_1 +\bar c\xi'_1
+\frac{1}{2}\Gamma_+\Big(- \sum_ix^i s_{\bar c,d}(e_i)d\xi_2
+ \sum_ix^i \bar cs_{\bar c,d}(e_i)\xi_2
+ \sum_ix^iB(e_i)\xi_2\Big)\\
&=
\big(\partial_-\xi'_1 +\bar c\xi'_1\big)
+\frac{1}{2}\Gamma_+\sum_i x^i 
\big( q_{\bar c,d}(e_i)+ B(e_i)\big)\xi_2\,.
\end{aligned}
\]
The vanishing of this term means that both linearly independent summands have to vanish separately for all $x^i$. 
Therefore, we end up with
\begin{gather}
\big(q_{\bar c,d}(v)+B(v)\big)\xi_2=0\,,\label{A}\\
\xi'_1(x^-)=\exp(-x^-\bar c)\xi_1^0\,,\label{B}
\end{gather}
with $\xi_2=\exp(-x^- d)\xi_2^0$ and (\ref{B}) following from $\partial_-\xi'_1 =-\bar c\xi'_1$. We recall that $\xi_1^0,\xi_2^0$ are constant spinors that obey $\sigma_+\xi_1^0=\sigma_-\xi_2^0=0$.
\begin{rem}\label{rem:parallel}
In terms of the local coordinates we see again that the space of parallel spinors is of dimension $\frac{1}{2}\dim S(W)$ and parametrized by $\xi_1^0$, generically. 
In case of maximal $\mathcal{K}_1$ we need the full freedom in the choice of $\xi_\alpha^0$ in (\ref{A}) and  (\ref{B}). In this case the vanishing of the bracket in (\ref{A}) is needed. This, of course, is the same as the vanishing of the sole remaining curvature term in (\ref{curv-MB}), i.e.\ (\ref{qcp}).
\end{rem}

\subsection{A family of eleven dimensional spaces}

From now on we are interested in non-flat connections and turn to 
dimension eleven. More precisely, we consider a connection that is given by\footnote{In \cite{klinker:CW} we showed that pairs of this type play an essential role, when we look for flat connections.} 
\begin{equation}\label{41}
\bar c:=(\alpha\Gamma_I+\beta\Gamma_J)\Gamma_K\,,\quad
     d:=(\alpha'\Gamma_I+\beta'\Gamma_J)\Gamma_K\,
\end{equation}
with $I,J,K\subset\{1,\ldots,9\}$ and $I\cap J\cap K=\emptyset$. We use projections 
\[
X_{IJ}^\pm:=\tfrac{1}{2}(\mathbbm{1}\pm \imath_{IJ}\Gamma_{IJ})
\]
with $\imath_{IJ}\in\{1,i\}$ such that $(\imath_{IJ}\Gamma_{IJ})^2=\mathbbm{1}$. In terms of $X^\pm_{IJ}$  we write
\begin{equation}\label{42}
q_{\bar c,d}(e_i)=\alpha_i^+\Gamma_iX_{IJ}^++\alpha_i^-\Gamma_i X_{IJ}^-
\end{equation}
for some linear combinations $\alpha_i^\pm\in\{\pm\alpha\pm\alpha'\pm\beta\pm\beta'\}$ where the specific arrangement of signs depend on whether $i\in I,J,K$, or $(I\cup J\cup K)^\complement$.   

We will further specify our connection and consider%
\footnote{
From the physics point of view the three-form -- or more precisely the combination $\Gamma_+s_{\bar c,d}(e_\mu)$ -- can be interpreted as the additional ingredient in the variation of the gravitino that comes from the curvature of the three-form potential in minimal eleven dimensional supergravity theories, see, for example, \cite{Cremm}.
From the mathematics point of view the Clifford elements of degree two or three in eleven dimensions are the only forms that are skew symmetric with respect to the charge conjugation. This is one essential condition for the general construction of super algebras from homogeneous connections, see \cite{klinker:habil}. Moreover, if we consider flat connections as in Proposition \ref{prop:flat} the three forms turn out to be the unique forms for which the geometric superalgebra turns into a supersymmetry algebra, see Examples \ref{ex:maximal}, \ref{ex:maximal-cont}, and references therein.
}
$|I|=|J|=|K|+1=2$ or -- without loss of generality -- $I=(12)$, $J=(34)$, and $K=(5)$, i.e.
\begin{equation}\label{43}
\bar c=(\alpha_+X_{1234}^++\alpha_-X_{1234}^-)\Gamma_{125}\,,\quad
d=(\alpha'_+X_{1234}^++\alpha'_-X_{1234}^-)\Gamma_{125}\,,
\end{equation}
with
\[
\alpha_\pm =\alpha\mp \beta,\quad \alpha_\pm'=\alpha'\mp\beta'\,.
\]
In this case $q_{\bar c,d}$ is given by
\begin{equation}\label{q}
\q_{\bar c,d}(e_i)=
\begin{cases}
(\alpha_--\alpha_+')^2\Gamma_iX^+_{1234}+(\alpha_+-\alpha_-')^2\Gamma_iX^-_{1234} & \text{for }i\in\{1,2\}\,,\\
(\alpha_-+\alpha_+')^2\Gamma_iX^+_{1234}+(\alpha_++\alpha_-')^2\Gamma_iX^-_{1234} & \text{for }i\in\{3,4\}\,,\\
(\alpha_+-\alpha_+')^2\Gamma_iX^+_{1234}+(\alpha_--\alpha_-')^2\Gamma_iX^-_{1234} & \text{for }i\in\{5\}\,,\\
(\alpha_++\alpha_+')^2\Gamma_iX^+_{1234}+(\alpha_-+\alpha_-')^2\Gamma_iX^-_{1234} & \text{for }i\in\{6,7,8,9\}\,.
\end{cases}
\end{equation}
\begin{rem}
In (\ref{A}) we have $X^{-}_{1234}\xi_2=\exp(-x^-d)X^{-}_{1234}\xi_2^0$ such that $X^{-}_{1234}\xi_2=0\iff X^{-}_{1234}\xi_2^0=0$. Therefore, condition (\ref{B}) is equivalent to 
\[
\big(q_{\bar c,d}(v)+B(v)\big)\xi_2^0=0\,.
\]
\end{rem}
Collecting the discussion yields the following proposition.
\begin{prop}\label{prop:11d-connection}
For the four parameter family of connections (\ref{41}) the space of parallel spinors is of dimension $\frac{3}{4}\dim S(W)$ and given by $S(V)\oplus X^+_{IJ}S(V)\subset S(W)$. The metric has at most four different eigenvalues determined by $\alpha_i^+$ from (\ref{42}). 

In the particular situation of (\ref{43}) this is true for the metric defined by\footnote{Below we will use 
$\lambda_1=\lambda_2=-i(\alpha_--\alpha_+')$, 
$\lambda_3=\lambda_4=-i(\alpha_-+\alpha_+')$, 
$\lambda_5=-i(\alpha_+-\alpha_+')$, and 
$\lambda_6=\ldots=\lambda_9=-i(\alpha_++\alpha_+')$.} 
\begin{equation}
B=-{\rm diag}\left(
(\alpha_--\alpha_+')^2\mathbbm{1}_2,
(\alpha_-+\alpha_+')^2\mathbbm{1}_2,
(\alpha_+-\alpha_+')^2\mathbbm{1}_1,
(\alpha_++\alpha_+')^2\mathbbm{1}_4  
\right)\,.\label{11d-metric-1}
\end{equation}
\end{prop}

\begin{rem}\label{rem:inv}
The Clifford elements $\bar c,d$ that define the connection are invariant with respect to $\mathfrak{so}_B(V)$, i.e.\ $[\bar c,A]=[d,A]=0$ for all $A\in\mathfrak{so}_B(V)\subset\cl(W)$. 
\end{rem}

The calculations in Section \ref{sec:Killing} show that the space of Killing vector fields $\mathcal{K}_0$ of the metric (\ref{11d-metric-1}) is spanned by 
\begin{equation}\label{K0}
\begin{aligned}
K_{(+)}& =-\partial_+\,,\quad K_{(-)}=-\partial_-\,, \\
K_{(i)}&= \cos(\lambda_ix^-)\partial_i + \lambda_i\sin(\lambda_1x^-)x^i\partial_+\,,\\
	%\begin{cases}
	%\cos(\lambda_1x^-)\partial_i + \lambda_1\sin(\lambda_1x^-)x^i\partial_+
	% & \text{for }i=1,2\,,\\
	%\cos(\lambda_2x^-)\partial_i + \lambda_2\sin(\lambda_2x^-)x^i\partial_+
	% & \text{for }i=3,4\,,\\
	%\cos(\lambda_3x^-)\partial_i + \lambda_3\sin(\lambda_3x^-)x^i\partial_+
	%& \text{for }i=5\,,\\
	%\cos(\lambda_4x^-)\partial_i + \lambda_4\sin(\lambda_4x^-)x^i\partial_+ 
	%& \text{for }6\leq i\leq 9	\,,
	%\end{cases}\\
K_{(i^*)} &= -\lambda_i\sin(\lambda_ix^-)\partial_i+\lambda_i^2 \cos(\lambda_ix^-)x^i\partial_+\,,\\
	%\begin{cases} 
	%-\lambda_1\sin(\lambda_1x^-)\partial_i+\lambda_1^2 \cos(\lambda_1x^-)x^i\partial_+
	% & \text{for }i=1,2\,,\\
	%-\lambda_2\sin(\lambda_2x^-)\partial_i+\lambda_2^2 \cos(\lambda_2x^-)x^i\partial_+
	% & \text{for }i=3,4\,,\\
	%-\lambda_3\sin(\lambda_3x^-)\partial_i+\lambda_3^2 \cos(\lambda_3x^-)x^i\partial_+
	%& \text{for }i=5\,,\\
	%-\lambda_4\sin(\lambda_4x^-)\partial_i+\lambda_4^2 \cos(\lambda_4x^-)x^i\partial_+
	%& \text{for }6\leq i\leq 9	\,,
	%\end{cases}\\
K_{(ij)} &= x^j\partial_i-x^i\partial_j\,,
\end{aligned}
\end{equation}
for $1\leq i\leq 9$ and $(ij)\in\{1,2\}^2\cup\{3,4\}^2\cup\{6,7,8,9\}^2$.

Analogously, for the space of spinors that are parallel with respect to the connection defined by the pair $(\bar c,d)$ according to (\ref{43}) the results of Section \ref{subsec:local}, namely (\ref{C})-(\ref{B}), yield
\begin{equation}\label{K1}
\mathcal{K}_1 = \left\{ \vec\xi\in\cancel S\,\left|\, 
\begin{array}{r@{\ }l}
\displaystyle \vec\xi=\vec\xi(\xi_1^0,\xi_2^0) = & \displaystyle  \exp(-x^-\bar c)\xi_1^0  \\[0.5ex]
& \displaystyle + \Big(1+\tfrac{1}{2}\sum_i\Gamma_+x^is_{\bar c,d}(e_i)\Big)\exp(-x^-d)\xi_2^0\,,
\\[1ex]
\multicolumn{2}{l}{\displaystyle \xi_1^0,\xi_2^0\text{ constant}\,,
\sigma_-\xi_2^0=\sigma_+\xi_1^0=0\,,
X^-_{1234}\xi_2^0=0}
\end{array}
\right.\right\}\,.
\end{equation}

\section{Geometric superalgebras}\label{sec:geoalg}

\subsection{Introduction}\label{sec:geoalg1}

In this section we will not give the definition of geometric superalgebras and geometric supersymmetry in general but consider the special situation from Proposition \ref{prop:11d-connection}. Nevertheless, we will shortly recall the idea. 

A manifold $M$ is said to admit a geometric superalgebra if there exists an extension of a Lie (sub)algebra  $\mathcal{K}_0(M)$ of the Killing vector fields on $M$ to a graded skew-symmetric superalgebra where the even part acts on the odd part in terms of derivations. 
The manifold is called to admit geometric supersymmetry if this extension is a Lie superalgebra. 
The odd part of this algebra is assumed to be purely geometric in the following sense: We consider a vector bundle $E$ over $M$ along with a connection $D$ such that $\mathcal{K}_0(M)$ acts on a subspace of the space of parallel sections 
\begin{equation}
\mathcal{K}_1(M)\subseteq\{s\in\cancel E\,|\,Ds=0\}.
\label{K1-general}
\end{equation}
This action then defines the even/odd-bracket of the (Lie) superalgebra. 
That means, there exits a map $\mathcal{L}\in{\rm Hom}(\mathcal{K}_0,{\rm End}(\mathcal{K}_1))$ with  $[\mathcal{L}_X,\mathcal{L}_Y]=\mathcal{L}_{[X,Y]}$ for all $X,Y\in\mathcal{K}_0(M)$. 
If the space $\mathcal{K}_1(M)$ in (\ref{K1-general}) is chosen maximal the extension is called non-restricted.

One non-trivial and important ingredient of such extension is the pairing $\{\cdot,\cdot\}:\mathcal{K}_1(M)\times\mathcal{K}_1(M)\to\mathcal{K}_0(M)$ that defines an algebra structure that is compatible with $\mathcal{L}$, i.e.\ $\big[X,\{s_1,s_2\}\big]= \big\{\mathcal{L}_Xs_1,s_2\big\}+\big\{s_1,\mathcal{L}_Xs_2\big\}$ for all $X\in\mathcal{K}_0(M)$, $s_1,s_2\in\mathcal{K}_1(M)$. 

The bundle $E$ usually is a spinor bundle over the base manifold $M$ and a geometric superalgebra or geometric supersymmetry is called irreducible if the spinor bundle is modeled on an irreducible Clifford module, say $S_0$, otherwise it is called reducible or $N$-extended. We will come back to this point in Section \ref{sec:7}, see Remark \ref{rem:N-ext}. 

In our special situation we consider the following data:

\begin{itemize}[leftmargin=4.5ex]
\item
The rank-$32$ bundle $S=S(M_B)$ of spinors over the eleven dimensional Cahen-Wallach space $M_B$ that is defined by (\ref{11d-metric-1}).
\item
The connection $D$ on $S$ defined by $(\bar c,d)$ according to (\ref{43}).
\item
The Lie algebra of Killing vector fields $\mathcal{K}_0$ and the space $\mathcal{K}_1\subset\cancel S$ of dimension $24$ that is given by the spinors parallel with respect to $D$, see (\ref{K0}) and (\ref{K1}).
\item
The basic ingredient to define the map $\mathcal{K}_1\times\mathcal{K}_1\to\mathcal{K}_{0}$ is a spin-invariant bilinear form on the space of spinors. Such bilinear forms has been fully classified, see\cite{AC,Che} for example. In dimension four such bilinear is used to transform a spinor that obeys the Dirac equation in such a way that it obeys the equation with opposite charge. Therefore, such bilinear form is frequently called charge conjugation and we will use this term, too.

The charge conjugation $C_W=C_V\otimes\sigma_2$ on the spinor space $S(W)=S(V)\otimes S_2$ of $W=V\oplus \RR^{1,1}$ is skew-symmetric and can be described in terms of the symmetric charge conjugation $C_V$ on the spinor space $S(V)$ of $V=\RR^9$ and the skew-symmetric charge conjugation $\sigma_2$ on $\RR^{1,1}$. 
This yields a symmetric map $S(W)\times S(W)\to W$  by $(\vec\xi,\vec\eta)\mapsto \big(C_W(\vec\xi,\Gamma_\mu\vec\eta)\big)_{\mu\in\{+,-,i\}}$. In terms of the two factors of $C_W$ this is
\begin{equation}\label{charge-flat}
\begin{gathered}
C_W(\vec\xi,\Gamma_-\vec\eta) = -\sqrt{2} iC_V(\xi_1,\eta_1)\,,\quad 
C_W(\vec\xi,\Gamma_+\vec\eta) = \sqrt{2} iC_V(\xi_2,\eta_2)\,,  \\
C_W(\vec\xi,\Gamma_i\vec\eta) = -iC_V(\xi_1,\gamma_i\eta_2)- i C_V(\eta_1,\gamma_i\xi_2)\,.
\end{gathered}
\end{equation}
The skew-symmetry of $C_W$ and the symmetry of the map $S(W)\otimes S(W)\to W$ implies that for a fixed Clifford element $a\in\cl(W)$ of degree $\ell$ the symmetry of $S(W)\times S(W)\to\Lambda^\ell W$ with $(\vec\xi,\vec\eta)\mapsto C_W(\vec\xi,a\,\vec\eta)=\Delta_\ell C_W(\vec\eta,a\,\vec\xi)$ is given by 
\begin{equation}\label{symmetry}
\Delta_\ell=-(-1)^{\frac{\ell(\ell+1)}{2}}\,.
\end{equation}
\item 
The charge conjugation on the standard fiber yields a spin-invariant bilinear form on the spinor bundle $S$ that we will denote by $C$, too. If we use the splitting of the bundle introduced at the beginning of Section \ref{subsec:local} and write $\vec\xi=\xi_1+\xi_2$ and $\vec\eta=\eta_1+\eta_2$ we have 
\[
C(\vec\xi,\vec\eta)=C(\xi_1,\eta_2)-C(\eta_1,\xi_2)\,.
\]
\item The action of the Killing vector fields $\mathcal{K}_0$ on the spinors is defined by the spinorial Lie derivative
\begin{equation}\label{spin-Lie}
\mathcal{L}:\mathcal{K}_0\times \cancel S\to \cancel S\,,\quad 
(K,\vec\xi)\mapsto \mathcal{L}_K\vec\xi:=\nabla_K\vec\xi -\Gamma(\nabla K)\vec\xi\,,
\end{equation}
see \cite{Kosmann:1971}. We emphasize the fact that this definition works properly only for Killing vector fields, because in this case $\nabla K$ is indeed skew symmetric. 

Nevertheless, it has to be checked, whether parallel spinors from $\mathcal{K}_1$ stay parallel after applying the Lie derivative, or, in other word, the connection is invariant under isometries.
\end{itemize}

\subsection{Even-Odd and Even-Even-Odd}\label{sec:eo}

\begin{prop}\label{prop:even-odd}
Consider $\mathcal{K}_0$ and $\mathcal{K}_1$ as defined in (\ref{K0}) and (\ref{K1}). 
The Lie derivative when restricted to $\mathcal{K}_0$ acts on $\mathcal{K}_1$.
\end{prop}

For the proof of Proposition \ref{prop:even-odd}, we need to know how the operators $\mathcal{L}_K=\nabla_K-\Gamma(\nabla K)$ for $K\in\mathcal{K}_0$ act on spinors from $\mathcal{K}_1$. We will do this for the basic elements $\{K_{(\pm)},K_{(i)},K_{(i^*)},K_{(ij)}\}$. 

We will make use of
\begin{align*}
\nabla_\mu(K_{(-)})_\nu 
	&= \partial_\mu(K_{(-)})_\nu + \sum_\kappa\Gamma_{\mu\kappa;\nu}(K_{(-)})^\kappa  
% &=\begin{cases}
%\sum_{j}\nu_{\ell j}x^j  & \text{for }(\mu\nu)=(-\ell)\,,\\
%-\sum_{j}\nu_{\ell j}x^j & \text{for }(\mu\nu)=(\ell-)\,,\\
%0 & \text{else}\,,
	 = \begin{cases}
		\lambda_\ell^2x^\ell & \text{for }(\mu\nu)=(-\ell)\,,\\
		-\lambda_\ell^2x^\ell & \text{for }(\mu\nu)=(\ell-)\,,\\
		0 & \text{else}\,,
	\end{cases}\\
\nabla_\mu(K_{(+)})_\nu 
	&= \partial_\mu(K_{(+)})_\nu + \sum_\kappa\Gamma_{\mu \kappa;\nu}(K_{(+)})^\kappa = 0\,,\\
\nabla_\mu(K_{(i)})_\nu 
	&= \partial_\mu(K_{(i)})_\nu + \Gamma_{\mu-;\nu}(K_{(i)})^- +\sum_j\Gamma_{\mu j;\nu}(K_{(i)})^j \\
	&= \partial_\mu(K_{(i)})_\nu + \sum_j\Gamma_{\mu j;\nu}(K_{(i)})^j 
%& =\begin{cases}
%\delta_{i\ell}\partial_-\alpha_\ell  & \text{for }(\mu\nu)=(-\ell)\,, \\
%\partial_-(\beta_ix^i) -\sum_{j\ell}B_{j\ell}x^\ell \alpha_j &\text{for } (\mu\nu)=(ii)\,,\\
%\partial_\ell (\beta_ix^i)  & (\mu\nu)=(\ell-)\,, \\
%0 &\text{else}
%\end{cases}\\
	 =\begin{cases}
		-\delta_{i\ell} \beta_i & \text{for }(\mu\nu)=(-\ell)\,, \\
		\delta_{i\ell} \beta_i  & \text{for }(\mu\nu)=(\ell-)\,, \\
		0 &\text{else}\,,
	\end{cases}\\
\nabla_\mu(K_{(i^*)})_\nu
	&=  \partial_\mu(K_{(i^*)})_\nu
	= \begin{cases}
		-\delta_{i\ell} \beta^*_i & \text{for }(\mu\nu)=(-\ell)\,, \\
		\delta_{i\ell} \beta^*_i  & \text{for }(\mu\nu)=(\ell-)\,, \\
		0 &\text{else}\,,
	\end{cases}\\
\intertext{and}
\nabla_\mu(K_{(ij)})_\nu 
	&=  \partial_\mu(K_{(ij)})_\nu +\sum_\ell \Gamma_{\mu\ell;\nu}(K_{(ij)})^\ell 
%\\&=\begin{cases} 
%	   \sum _\ell\Gamma_{-\ell;-}K_{(ij)}^\ell & \text{for }(\mu\nu)=(--) \\
%	   \delta_{kj}\delta_{i\ell}-\delta_{ik}\delta_{j\ell}  & \text{for }(\mu\nu)=(k\ell)\\
%	\end{cases}\\
%&=\begin{cases} 
%	   \sum_{k}\delta_{j k}\lambda_j^2 x^kx^i
%	   -\sum_{k}\delta_{i k}\lambda_i^2 x^kx^j & \text{for }(\mu\nu)=(--) \\
%	   \delta_{kj}\delta_{i\ell}-\delta_{ik}\delta_{j\ell}  & \text{for }(\mu\nu)=(k\ell)\\
%	\end{cases}\\&
=\begin{cases}    
	   \delta_{kj}\delta_{i\ell}-\delta_{ik}\delta_{j\ell}  & \text{for }(\mu\nu)=(k\ell)\,,\\
	   0 &\text{else}\,.
	\end{cases}
\end{align*}
We will also use 
\begin{equation*}
\begin{aligned}
s_{\bar c,d}(e_i)&=
%\begin{cases}
%\big((\alpha_+-\alpha_-')X^+_{1234}+(\alpha_--\alpha_+')X^-_{1234}\big)\Gamma_{125}\Gamma_i& i=1,2\\
%\big((\alpha_++\alpha_-')X^+_{1234}+(\alpha_-+\alpha_+')X^-_{1234}\big)\Gamma_{125}\Gamma_i& i=3,4\\
%\big((\alpha_+-\alpha_+')X^+_{1234}+(\alpha_--\alpha_-')X^-_{1234}\big)\Gamma_{125}\Gamma_i& i=5 \\
%\big((\alpha_++\alpha_+')X^+_{1234}+(\alpha_-+\alpha_-')X^-_{1234}\big)\Gamma_{125}\Gamma_i& i=6\ldots 9
%\end{cases}\\
%&=
\begin{cases}
 (\alpha_--\alpha_+')\Gamma_{125}\Gamma_iX^+_{1234}
+(\alpha_+-\alpha_-')\Gamma_{125}\Gamma_iX^-_{1234}& i=1,2\\
 (\alpha_-+\alpha_+')\Gamma_{125}\Gamma_iX^+_{1234}
+(\alpha_++\alpha_-')\Gamma_{125}\Gamma_iX^-_{1234}& i=3,4\\\
 (\alpha_+-\alpha_+')\Gamma_{125}\Gamma_iX^+_{1234}
+(\alpha_--\alpha_-')\Gamma_{125}\Gamma_iX^-_{1234}& i=5\ \\
 (\alpha_++\alpha_+')\Gamma_{125}\Gamma_iX^+_{1234}
+(\alpha_-+\alpha_-')\Gamma_{125}\Gamma_iX^-_{1234}& i=6,\ldots,9
\end{cases}
\end{aligned}
\end{equation*}
which for $\vec\xi(\xi_1^0,\xi_2^0)$ according to (\ref{K1}) becomes
\[
s_{\bar c,d}(e_i)\vec\xi = i\lambda_i\Gamma_{125}\Gamma_iX^+_{1234}\vec\xi\,.
\]
Furthermore, we will need
\begin{align*}
&\begin{aligned}
\exp(-x^-d)=\ &\big(\cos(i\alpha_+'x^-)-i\sin(i\alpha_+'x^-)\Gamma_{125}\big) X^+_{1234} \\
			 & +\big(\cos(i\alpha_-'x^-)-i\sin(i\alpha_-'x^-)\Gamma_{125}\big) X^-_{1234}\,,
\end{aligned}\\
&\begin{aligned}
\exp(-x^-c)=\ & \big(\cos(i\alpha_+x^-)-i\sin(i\alpha_+x^-)\Gamma_{125}\big) X^+_{1234}\\
			  &+\big(\cos(i\alpha_-x^-)-i\sin(i\alpha_-x^-)\Gamma_{125}\big) X^-_{1234}\,.
\end{aligned}			  
\end{align*}
A careful calculation yields
\begin{align*}
\mathcal{L}_{K_{(+)}}\vec\xi
 &=\nabla_+\vec\xi=0\,,\\
\mathcal{L}_{K_{(-)}}\vec\xi
 &=\nabla_-\vec\xi - \frac{1}{2}\sum_j\lambda_j^2x^j\Gamma_+\Gamma_j\vec\xi \\
 & = -\bar c\xi_1- d \xi_2 - \frac{1}{2}\sum_j\lambda_j^2x^j\Gamma_+\vec\xi \\
 & = -\exp(-x^-\bar c) \bar c\xi_1^0- \exp(-x^-d) d \xi_2^0  
 -\frac{1}{2}\Gamma_+ \sum_jx^j\big( \bar c s_{\bar c,d}(e_j)+\lambda_j^2\Gamma_j\big)\xi_2\\
 &= -\exp(-x^-\bar c) \bar c\xi_1^0- \exp(-x^-d) d \xi_2^0  
 -\frac{1}{2}\Gamma_+ \sum_jx^j s_{\bar c,d}(e_j) d\xi_2 \\
 & = -\vec\xi(\bar c \xi_1^0,d\xi_2^0)\,,
\end{align*}
as well as 
\begin{align*}
&\begin{aligned}
\mathcal{L}_{K_{(i)}}\vec\xi
	&=\alpha_i(x^-)\nabla_i\vec\xi +\beta_i(x^-)x^i\nabla_+\vec\xi 
		+\frac{1}{2}\Gamma_+\beta_i\Gamma_i\vec\xi\\
	&= \frac{1}{2}\Gamma_+ \big( \alpha_i(x^-)s_{\bar c,d}(e_i)+\beta_i(x^-)\Gamma_i\big)\xi_2\\
 	&= \frac{i\lambda_i}{2}\Gamma_+ \big( \cos(\lambda_ix^-)\Gamma_{125}\Gamma_i
 		-i\sin(\lambda_ix^-)\Gamma_i\big)\xi_2\,.
\end{aligned}
\end{align*}
On the one hand -- by recalling $X^+_{1234}\xi_2=\xi_2$ -- we have
\begin{align*}
\mathcal{L}_{K_{(i)}}\vec\xi
	=\ &\frac{i\lambda_i}{2}\Gamma_+ \big( \cos(\lambda_ix^-) \Gamma_{125}\Gamma_i
			-i\sin(\lambda_ix^-)\Gamma_i\big)
			\big(\cos(i\alpha_+'x^-)-i\sin(i\alpha_+'x^-)\Gamma_{125}\big)\xi^0_2\\
%	=\ & \frac{i\lambda_i}{2}\Gamma_+ \Big( 
%		    \big(\cos(\lambda_ix^-)\cos(i\alpha_+'x^-)
%		     	-\epsilon_i\sin(i\alpha_+'x^-)\sin(\lambda_ix^-)\big)   \Gamma_{125}\Gamma_i\\
%	&	    +\big( -i\epsilon_i\cos(\lambda_ix^-)\sin(i\alpha_+'x^-)
%		    	-i\sin(\lambda_ix^-)\cos(i\alpha_+'x^-) \big)  \Gamma_i \Big)\\
	=\ & \begin{cases}
			\frac{i\lambda_i}{2}\Gamma_+ \big( 
			\cos( (\lambda_i+i\alpha_+')x^-)\Gamma_{125}\Gamma_i
			-i\sin((\lambda_i+i\alpha_+')x^-)\Gamma_i
		\big)\xi_2^0 &\text{for } i=1,2,5\\[1.5ex]
\frac{i\lambda_i}{2}\Gamma_+ \big( 
			\cos( (\lambda_i-i\alpha_+')x^-)\Gamma_{125}\Gamma_i
			-i\sin((\lambda_i-i\alpha_+')x^-)\Gamma_i
		\big)\xi_2^0 &\text{else } \\		
	\end{cases}	\\	
	=\ & \begin{cases}
\frac{i\lambda_i}{2}\Gamma_+ \big( 
			\cos( i\alpha_-x^-)\Gamma_{125}\Gamma_i+i\sin(i\alpha_-x^-)\Gamma_i
			\big)\xi_2^0 &\text{for } i=1,\ldots,4\\[1.5ex]		
\frac{i\lambda_i}{2}\Gamma_+ \big( 
			\cos( i\alpha_+x^-)\Gamma_{125}\Gamma_i+i\sin(i\alpha_+x^-)\Gamma_i
		\big)\xi_2^0 &\text{for } i=5,\ldots,9		
	\end{cases}	
\end{align*}
and, on the other hand we have
\begin{align*}
\vec\xi\big(s_{\bar c,d}(e_i)\xi_2^0,0\big) 
=\ &\exp(-cx^-)s_{\bar c,d}(e_i)\xi_2^0 \\
=\ &i\lambda_i \big(\cos(i\alpha_+x^-)-i\sin(i\alpha_+x^-)\Gamma_{125}\big) X^+_{1234}   \Gamma_i\Gamma_{125}\xi_2^0 \\
& +i\lambda_i\big(\cos(i\alpha_-x^-)-i\sin(i\alpha_-x^-)\Gamma_{125}\big) X^-_{1234}\Gamma_i\Gamma_{125}\xi_2^0 \\
=\ &\begin{cases}
i\lambda_i \big(\cos(i\alpha_-x^-)-i\sin(i\alpha_-x^-)\Gamma_{125}\big)\Gamma_{125}\Gamma_i\xi_2^0 &\text{for } i=1,\ldots,4 \\[1,5ex]
i\lambda_i\big(\cos(i\alpha_+x^-)-i\sin(i\alpha_+x^-)\Gamma_{125}\big) \Gamma_{125}\Gamma_i\xi_2^0 &\text{for } i=5,\ldots,6 
\end{cases}
\end{align*}
such that 
\begin{align*}
\mathcal{L}_{K_{(i)}}\vec\xi
	=\ & \frac{1}{2}\vec\xi\big(\Gamma_+s_{\bar c,d}(e_i)\xi_2^0,0\big)\,.
\end{align*}
Doing analogous calculations for 
\begin{align*}
&\begin{aligned}
\mathcal{L}_{K_{(i^*)}}\vec\xi
	&= \frac{1}{2}\Gamma_+ \big( \alpha^*_i(x^-)s_{\bar c,d}(e_i)+\beta^*_i(x^-)\Gamma_i\big)\xi_2 \\
 	&= \frac{\lambda_i^2}{2}\Gamma_+ \big(\cos(\lambda_ix^-)\Gamma_i 
 			- i\sin(\lambda_ix^-)\Gamma_{125}\Gamma_i\big)\xi_2
\end{aligned}
\end{align*}
we get 
\[
\mathcal{L}_{K_{(i^*)}}\vec\xi=-\frac{\lambda_i^2}{2}\Gamma_+\vec\xi\big(\Gamma_i\xi_2^0,0\big)\,.
\]
In addition we obtain
\[
\mathcal{L}_{K_{(ij)}}\vec\xi=\frac{1}{2}\vec\xi\big(\Gamma_{ij}\xi_1^0,\Gamma_{ij}\xi_2^0\big)
\]
where we have to take into account $(ij)\in\{1,2\}^2\cup\{3,4\}^2\cup\{6,\ldots,9\}^2$.
We collect the result in the following remark.
\begin{rem}\label{rem:liederivative}
The Lie derivatives according to Proposition \ref{prop:even-odd} are explicitly given by
\begin{align*}
\mathcal{L}_{K_{(+)}}\vec\xi\big(\xi_1^0,\xi_2^0\big) &= 0\,,\\
\mathcal{L}_{K_{(-)}}\vec\xi\big(\xi_1^0,\xi_2^0\big)	&= -\vec\xi\big(\bar c\xi_1^0,d \xi_2^0\big)\,,\\
\mathcal{L}_{K_{(i)}}\vec\xi\big(\xi_1^0,\xi_2^0\big)	
		&= -\vec\xi\big(-\tfrac{1}{2}\Gamma_+s_{\bar c,d}(e_i)\xi_2^0,0\big)\,,\\
\mathcal{L}_{K_{(i^*)}}\vec\xi\big(\xi_1^0,\xi_2^0\big)
		&= -\vec\xi\big(\tfrac{1}{2}\Gamma_+B(e_i)\xi_2^0,0\big)\,,\\
\mathcal{L}_{K_{(ij)}}\vec\xi\big(\xi_1^0,\xi_2^0\big)
		&= -\vec\xi\big(-\tfrac{1}{2}\Gamma_{ij}\xi_1^0,-\tfrac{1}{2}\Gamma_{ij}\xi_2^0\big)\,.
\end{align*}
For $\mu\in\{\pm,i,i^*,ij\}$ we can rewrite this as 
\begin{equation}\label{50}
\mathcal{L}_{K_{(\mu)}}\vec\xi\big(\xi_1^0,\xi_2^0\big) = -\vec\xi\big(\rho(e_\mu)(\xi_1^0,\xi_2^0)\big)
\end{equation}
with $\rho$ according to Proposition \ref{prop:maps} and Remark \ref{rem:soBV} as well as $e_{i^*}:=e_i^*$ and $\mathfrak{so}_B(V)={\rm span}\{e_{ij}\}$ with $\rho(e_{ij})=-\frac{1}{2}\Gamma_{ij}$.
\end{rem}

\subsection{Odd-Odd and Even-Odd-Odd}

\begin{defn}\label{def:odd-odd}
We use the charge conjugation $C$ on $S$ to define a symmetric map
\[
\mathcal{K}_1\otimes\mathcal{K}_1\to   \mathcal{K}_0\,.
\]
Motivated by (\ref{charge-flat}) we consider the projection
\[
\mathcal{K}_1\otimes\mathcal{K}_1\to {\rm span}\{K_{(+)},K_{(-)},K_{(i)}\}\subset  \mathcal{K}_0\,.
\]
given by 
\begin{align}
\{\vec\xi,\vec\eta\}^W &= \{\vec\xi,\vec\eta\}^+ K_{(+)}+\{\vec\xi,\vec\eta\}^- K_{(-)}+ \sum_i \{\vec\xi,\vec\eta\}^i K_{(i)}
\end{align}
with
\begin{equation}\label{proj:W}\begin{gathered}
\{\vec\xi,\vec\eta\}^- = C(\xi^0_2,\Gamma_+\eta_2^0)\,,\quad
\{\vec\xi,\vec\eta\}^+ = C(\xi^0_1,\Gamma_-\eta_1^0)\,,\\
\{\vec\xi,\vec\eta\}^i = C(\xi^0_1,\Gamma_i\eta_2^0)+C(\eta^0_1,\Gamma_i\xi_2^0) \,.	
\end{gathered}\end{equation}
We complete this projection to $\{\cdot,\cdot\}$ by the two maps
\begin{align}
\mathcal{K}_1\otimes\mathcal{K}_1 \to\ & {\rm span}\{K_{(i^*)}\}\subset  \mathcal{K}_0\,,\quad
\{\vec\xi,\vec\eta\}^{*}  = \sum_i \{\vec\xi,\vec\eta\}^{i^*} K_{(i^*)}\,,
\\\intertext{and}
\mathcal{K}_1\otimes\mathcal{K}_1 \to\ & {\rm span}\{K_{(ij)}\} \subset  \mathcal{K}_0\,,\quad
\{\vec\xi,\vec\eta\}^{\mathfrak{so}} = \frac{1}{2}\sum_{ij} \{\vec\xi,\vec\eta\}^{ij} K_{(ij)}\,.
\end{align}
The coefficients therein are defined by
\begin{equation}\label{proj:V*}
\begin{aligned}
\{\vec\xi,\vec\eta\}^{i^*} =\ &
 		C\big(\xi_1^0, s_{\bar c,d}(B^{-1}(e_i))\eta_2^0\big)
 		-C\big(s_{\bar c,d}(B^{-1}(e_i)) \xi_2^0, \eta_1^0\big) \\
 =\ &  	 \frac{i}{\lambda_i}C\big( \xi_1^0,  \Gamma_{125}\Gamma_i\eta_2^0\big)
 		+\frac{i}{\lambda_i}C\big(\eta_1^0,  \Gamma_{125}\Gamma_i \xi_2^0\big) \\
\end{aligned}
\end{equation}
and
\begin{equation}\label{proj:so}
\begin{aligned}
\{\vec\xi,\vec\eta\}^{ij} =\ &
-\frac{1}{2} C\big(\xi^0_2, \Gamma_+ \big(s_{d,\bar c}(e_j)\Gamma_i+\Gamma_{i}s_{\bar c,d}(e_j)\big)\eta^0_2\big) \\
=\ & 
\begin{cases}
i\lambda_1  C\big(\xi^0_2,\Gamma_+\Gamma_{5}\eta_2^0\big) & \text{for }(ij)=(12)\\
i\lambda_3  C\big(\xi^0_2,\Gamma_+\Gamma_{12345}\eta_2^0\big) & \text{for }(ij)=(34)\\
i\lambda_6 C\big(\xi^0_2,\Gamma_+\Gamma_{125ij}\eta_2^0\big) & \text{for }(ij)\in\{6,\ldots,9\}^2
\end{cases}
\\
=:\ &i\epsilon_j\lambda_j C\big(\xi^0_2,\Gamma_+\Gamma_{125}\Gamma_{ij}\eta_2^0\big)\ \text{ for }(ij)\in\{1,2\}^2\cup\{3,4\}^2\cup\{6,\ldots,9\}^2
\end{aligned}
\end{equation}
In (\ref{proj:W}),(\ref{proj:V*}) and (\ref{proj:so}) the spinors $\xi_1^0,\xi_2^0$ and $\eta_1^0,\eta_2^0$ are the constant spinors that define the parallel spinors $\vec\xi$ and $\vec\eta$, see Remark \ref{rem:parallel} and the calculations before. 
\end{defn}

\begin{rem}
The construction in Definition \ref{def:odd-odd} can be made more general such that it provides a superalgebra for a wide class of connections according to Proposition \ref{prop:maps}. More precisely, it can be shown that this is the only possible algebra structure. This is used in \cite{klinker:habil} to start a systematic classification of supersymmetric extensions of Cahen-Wallach spaces that in particular covers the  examples found in the literature, see, for example,  \cite{CheKo84,Pope:2002,FigPapado1,Gauntlett:2002,Fig1,Fig03,MeFig04,Hustler,Meessen2002}. A first attempt of such systematic treatment has been started in \cite{santi1} with odd-odd bracket $\{\cdot,\cdot\}=\{\cdot,\cdot\}^W$. This definition turned out to be too restrictive for allowing a non-trivial superalgebra.
\end{rem}

\begin{prop}\label{prop:e-o-o}
For any $K\in\mathcal{K}_0$ and $\xi\in\mathcal{K}_1$ 
\[
\big[K,\{\vec\xi,\vec\xi\}\big]=2\{\mathcal{L}_K\vec\xi,\vec\xi\}\,.
\]
\end{prop}

The statement of Proposition \ref{prop:e-o-o} is clear for $K=K_{(+)}$ such that we may restrict ourselves to $K\in\{K_{(-)},K_{(i)},K_{(i^*)},K_{(ij)}\}$. The remaining proof needs the symmetry of the charge conjugation as stated in (\ref{symmetry}). For $\vec\xi=\vec\xi(\xi_1^0,0)$ we have 
\[
[K,\{\vec\xi,\vec\xi\}]=C(\xi_1^0,\Gamma_-\xi_1^0)[K,K_{(+)}]=0\
\text{ and }\
\{\mathcal{L}_K\vec\xi,\vec\xi\}=0\,.
\]
The last equation is only non obvious for $K=K_{(-)}, K_{(ij)}$ and is then due to
\begin{align*}
C(\bar c\xi_1^0,\Gamma_-\xi_1^0 )  
	&= \Delta_0\Delta_3 C(\xi_1^0,\bar c\Gamma_-\xi_1^0) 
	=  C(\xi_1^0,\bar c\Gamma_-\xi_1^0) \\
	&= \Delta_4 C(\xi_1^0,\bar c\Gamma_-\xi_1^0) 
	= -C(\xi_1^0,\bar c\Gamma_-\xi_1^0)\,, \\
C(\Gamma_{ij}\xi_1^0,\Gamma_-\xi_1^0 )
	&= \Delta_0\Delta_2 C(\xi_1^0,\Gamma_{ij}\Gamma_-\xi_1^0)
	= - C(\xi_1^0,\Gamma_{ij}\Gamma_-\xi_1^0)\\
	&= -\Delta_3 C(\xi_1^0,\Gamma_{ij}\Gamma_-\xi_1^0)
	= C(\xi_1^0,\Gamma_{ij}\Gamma_-\xi_1^0)\,.
\end{align*}
If we consider $\vec\xi=\vec\xi(0,\xi_2^0)$ we get
\begin{align*}
[K_{(-)},\{\vec\xi,\vec\xi\}] \ &= \{\vec\xi,\vec\xi\}^-[K_{(-)},K_{(-)}]  
+\frac{1}{2}\sum_{ij}\{\vec\xi,\vec\xi\}^{ij} [K_{(-)},K_{(ij)}] =0\,,
\\ 
\{ \mathcal{L}_{K_{(-)}}\vec\xi,\vec\xi\}
=\ & \{ \vec\xi(0,d\xi_2^0),\vec\xi(0,\xi_2^0) \}\\
=\ &  C(\xi^0_2,\Gamma_+d\xi_2^0) K_{(-)} 
		+ \frac{i\lambda_6}{2} \sum_{ij=6}^9 C(\xi^0_2,\Gamma_+\Gamma_{125ij}d\xi_2^0)K_{(ij)}\\
&	 	+i\lambda_1 C(\xi^0_2,\Gamma_+\Gamma_5d\xi_2^0)K_{(12)}
	 	+i\lambda_3 C(\xi^0_2,\Gamma_+\Gamma_{12345}d\xi_2^0)K_{(34)}\\	 
=\ &  0	 \,,
\end{align*}
because for $(ij)\in \{3,4\}^2\cup\{6,\ldots,9\}^2$ we have
\begin{align*}
C(\xi^0_2,\Gamma_+\Gamma_{125ij}d\xi_2^0) 
	&= \Delta_6 C(d\xi^0_2,\Gamma_+\Gamma_{125ij}\xi_2^0)
	= \Delta_6\Delta_0\Delta_3 C(\xi^0_2,d\Gamma_+\Gamma_{125ij}\xi_2^0)\\
	&= -\Delta_6\Delta_0\Delta_3 C(\xi^0_2,\Gamma_+\Gamma_{125ij}d\xi_2^0)
	= -C(\xi^0_2,\Gamma_+\Gamma_{125ij}d\xi_2^0)\,,\\
C(\xi^0_2,\Gamma_+\Gamma_{5}d\xi_2^0) 
	&= \Delta_2 C(d\xi^0_2,\Gamma_+\Gamma_{5}\xi_2^0)
	= \Delta_2\Delta_0\Delta_3 C(\xi^0_2,d\Gamma_+\Gamma_{5}\xi_2^0)\\
	&= -\Delta_2\Delta_0\Delta_3 C(\xi^0_2,\Gamma_+\Gamma_{5}d\xi_2^0)
	= - C(\xi^0_2,\Gamma_+\Gamma_{5}d\xi_2^0)\,,\\
C(\xi^0_2,\Gamma_+d\xi_2^0)	
	&= \Delta_1 C(d\xi^0_2,\Gamma_+\xi_2^0)
	= \Delta_1\Delta_0\Delta_3 C(\xi^0_2,d\Gamma_+\xi_2^0)\\
	&= -\Delta_1\Delta_0\Delta_3 C(\xi^0_2,\Gamma_+d\xi_2^0)
	= - C(\xi^0_2,\Gamma_+\Gamma_{5}d\xi_2^0)\,.
\end{align*}
For $K=K_{(k\ell)}$ we get 
\begin{align*}
[K_{(k\ell)},\{\vec\xi,\vec\xi\}]
=\ & \{\vec\xi,\vec\xi\}^-[K_{(k\ell)},K_{(-)}]  
	+\frac{1}{2}\sum_{ij}\{\vec\xi,\vec\xi\}^{ij} [K_{(k\ell)},K_{(ij)}] \\	
=\ & - i \sum_{ij}\epsilon_j\lambda_j C(\xi^0_2,\Gamma_+\Gamma_{125}\Gamma_{ij}\xi_2^0)\delta_{i[k}K_{(\ell] j)}
		\\
=\ &\begin{cases}
0 & \text{for }(k\ell)=(12),(34)\\
 i\lambda_6 \sum_{ij=6}^9 C(\xi^0_2,\Gamma_+\Gamma_{125}\Gamma_{j[k}\xi_2^0) K_{(\ell] j)}& \text{for }(k\ell)\in\{6,\ldots,9\}^2
\end{cases}\\
2\{\mathcal{L}_{K_{(k\ell)}}\vec\xi,\vec\xi\}
=\ & \{\vec\xi(0,\Gamma_{k\ell}\xi_2^0),\vec\xi(0,\xi_2^0)\} \\
=\ & \underbrace{C(\xi^0_2,\Gamma_+\Gamma_{k\ell}\xi_2^0)}_{=0,\ \Delta_3=-1} K_{(-)} 
	+ \frac{i\lambda_6}{2} \sum_{ij=6}^9 C(\xi^0_2,\Gamma_+\Gamma_{125ij}\Gamma_{k\ell}\xi_2^0)
			K_{(ij)}\\
&	 +i\lambda_1 \underbrace{C(\xi^0_2,\Gamma_+\Gamma_5\Gamma_{k\ell}\xi_2^0)}_{=0,\ \Delta_4=-1}K_{(12)}
	 +i\lambda_3 \underbrace{C(\xi^0_2,\Gamma_+\Gamma_{12345}\Gamma_{k\ell}\xi_2^0)}_{=0,\ \Delta_4=\Delta_8=-1}K_{(34)}\\
=\ & \frac{i\lambda_6}{2} \sum_{ij=6}^9 C(\xi^0_2,\Gamma_+\Gamma_{125}\Gamma_{ij}\Gamma_{k\ell}\xi_2^0)
			K_{(ij)}\\	
=\ & 
\begin{cases}
0 & \text{for }(k\ell)=(12),(34)\\
 i\lambda_6 \sum_{ij=6}^9 C(\xi^0_2,\Gamma_+\Gamma_{125}\Gamma_{j[k}\xi_2^0)
			K_{(\ell]j)} & \text{for }(k\ell)\in\{6,\ldots,9\}^2		
\end{cases}				
\end{align*}
For $K=K_{(k)}$ we get
\begin{align*}
[K_{(k)},\{\vec\xi,\vec\xi\}] 
=\ & \{\vec\xi,\vec\xi\}^-[K_{(k)},K_{(-)}]  + \frac{1}{2}\sum_{ij}\{\vec\xi,\vec\xi\}^{ij}[K_{(k)},K_{(ij)}]\\
=\ & C(\xi_2^0,\Gamma_+\xi_2^+) K_{(k^*)}  +\frac{1}{2}\sum_{ij}\{\vec\xi,\vec\xi\}_{ij}(\delta_{kj}K_{(i)}-\delta_{ki}K_{(j)})\\
=\ & C(\xi_2^0,\Gamma_+\xi_2^+) K_{(k^*)}  
  + \{\vec\xi,\vec\xi\}_{12}(\delta_{k2}K_{(1)}-\delta_{k1}K_{(2)})\\
& + \{\vec\xi,\vec\xi\}_{34}(\delta_{k4}K_{(3)}-\delta_{k3}K_{(4)})
  + \sum_{j=6}^9\{\vec\xi,\vec\xi\}_{jk}K_{(j)}\,,\\
2\{\mathcal{L}_{K_{(k)}}\vec\xi,\vec\xi\}
=\ & \{\vec\xi(\Gamma_+s_{\bar c,d}(e_k)\xi_2^0,0),\vec\xi(0,\xi_2^0)\} \\
=\ & \sum_j \frac{i}{\lambda_j} C(\Gamma_+s_{\bar c,d}(e_k)\xi_2^0,\Gamma_{125}\Gamma_j\xi_2^0) K_{(j^*)}\\
&	+ \sum_j C(\Gamma_+s_{\bar c,d}(e_k)\xi_2^0, \Gamma_j\xi_2^0)K_{(j)} \\
%=\ & - \sum_j \frac{\lambda_k}{\lambda_j} C(\Gamma_+\Gamma_{125}\Gamma_k X^+_{1234}\xi_2^0,
%	\Gamma_{125}\Gamma_j\xi_2^0) K_{(j^*)} \\
%&  + i\lambda_k \sum_j C(\Gamma_+\Gamma_{125}\Gamma_k X^+_{1234}\xi_2^0, \Gamma_j\xi_2^0)K_{(j)} \\
=\ & - \sum_j \frac{\lambda_k}{\lambda_j} C( \xi_2^0,\Gamma_+X^+_{1234} \Gamma_k\Gamma_j X^+_{1234}\xi_2^0)
 K_{(j^*)} \\
&  + i\lambda_k \sum_j C(\xi_2^0,\Gamma_+ X^+_{1234} \Gamma_k \Gamma_{125}\Gamma_j\xi_2^0)K_{(j)} \\
=\ &  C( \xi_2^0,\Gamma_+ \xi_2^0) K_{(j^*)} \\
&  + i\lambda_1 C(\xi_2^0,\Gamma_+  \Gamma_{5} \xi_2^0)(\delta_{k2}K_{(1)}-\delta_{k1}K_{(2)}) \\
&  + i\lambda_3 C(\xi_2^0,\Gamma_+  \Gamma_{12345}\xi_2^0)(\delta_{k4}K_{(3)}-\delta_{k3}K_{(4)}) \\
&  + i\lambda_6 \sum_{j=6}^9 C(\xi_2^0,\Gamma_+ \Gamma_{125jk}\xi_2^0)K_{(j)} \,.
\end{align*}
In the last step all other summands vanish because they are skew-symmetric, e.g.\ 
\[
\lambda_3 C(\xi_2^0,\Gamma_+ X^+_{1234} \Gamma_4\Gamma_{125}\Gamma_5\xi_2^0)K_{(5)}
= -\lambda_3 C(\xi_2^0,\Gamma_+ \Gamma_{124} \xi_2^0)K_{(5)} \overset{\Delta_4=-1}{=} 0\,, 
\] 
or because $X^-_{1234}\xi_2^0=0$, e.g.\ 
\[
\lambda_1 C(\xi_2^0,\Gamma_+  X^+_{1234}\Gamma_5\Gamma_{125}\Gamma_1 \xi_2^0)K_{(1)} 
= -\lambda_1 C(\xi_2^0,\Gamma_+  \Gamma_{2} X^-_{1234}\xi_2^0)K_{(1)} =0\,.
\]
By similar arguments we get that for $K=K_{(k^*)}$ the following terms coincide:
\begin{align*}
[K_{(k^*)},\{\vec\xi,\vec\xi\}] 
=\ & \{\vec\xi,\vec\xi\}^-[K_{(k^*)},K_{(-)}]  + \frac{1}{2}\sum_{ij}\{\vec\xi,\vec\xi\}^{ij}[K_{(k^*)},K_{(ij)}]\,,\\
2\{\mathcal{L}_{K_{(k^*)}}\vec\xi,\vec\xi\}
=\ & -\lambda^2_k\{\vec\xi(\Gamma_+\Gamma_k\xi_2^0,0),\vec\xi(0,\xi_2^0)\}\,.
\end{align*}

We collect the results of Section \ref{sec:geoalg} in the following statement.

\begin{thm}\label{thm:geoalg}
The indecomposable Cahen-Wallach space $M_B$ along with the connection $D$ according to Proposition \ref{prop:11d-connection} define a non-restricted geometric superalgebra. 
\begin{itemize}
\item The even and odd parts, $\mathcal{K}_0$ and $\mathcal{K}_1$, of the underlying graded vector space are given by (\ref{K0}) and (\ref{K1}), respectively.
\item The product structure is given by the usual commutator on $\mathcal{K}_0$ and completed by the even-odd bracket defined by the Lie derivative, see Remark \ref{rem:liederivative}, and the odd-odd bracket according to (\ref{proj:W}), (\ref{proj:V*}), and (\ref{proj:so}), from Definition \ref{def:odd-odd}.
\end{itemize}
\end{thm}

We end this section with a short comment on the question if the odd part of the geometric superalgbra is minimal in a certain sense.

\begin{rem}\label{rem:singular}
Because $K_{(i^*)}$ acts by $\mathcal{L}_{K_{(i^*)}}\vec\xi(0,\xi_2^0)=\lambda_i\vec\xi(\Gamma_+\Gamma_i\xi_2^0,0)$ we have 
\begin{equation*}
\mathcal{L}_{K_{(i^*)}}\vec\xi(0,\xi_2^0)\in 
\begin{cases}
X^+_{1234}S&\text{for }i=5,\ldots,9\,,\\ 
X^-_{1234}S&\text{for }i=1,\ldots,4\,.
\end{cases}
\end{equation*}
Therefore, a reduction of $\mathcal{K}_1$ in Theorem \ref{thm:geoalg} is only possible if some of the eigenvalues of $B$ vanish, i.e.\ if $M_B$ is decomposable.\footnote{As reduction we consider only those for which the resulting algebra is nontrivial, i.e. $\Gamma_+\mathcal{K}_1\neq 0$.} 

There are special configurations of parameters in the metric (\ref{11d-metric-1}) that yield decomposable spaces, namely $\alpha_+=\pm\alpha_+'$ and $\alpha_-=\pm\alpha_+'$.
In fact a reduction is only possible if $\alpha_-=\alpha_+'=0$ or $\alpha_+=\alpha_+'=0$.
\end{rem}

\section{Geometric Supersymmetry}\label{sec:susy}

\subsection{Odd-Odd-Odd}

In this section we will show that for a special set of parameters the geometric superalgebra from Theorem \ref{thm:geoalg} in fact defines geometric supersymmetry.

We recall the fact that a superalgebra is a Lie superalgebra if the graded Jacobi identity is fulfilled, i.e.\ for all elements $x,y,z$ we have
\[
(-1)^{|x||z|}[x,[y,z]+(-1)^{|y||z|}[z,[x,y]]+(-1)^{|x||y|}[y,[z,x]]=0
\]
where $|\cdot|$ denotes the $\ZZ_2$-degree. 

If a Cahen-Wallach space $M_B$, or, more precisely, $\mathcal{K}_0\oplus\mathcal{K}_1$, defines a geometric superalgebra the only obstruction to geometric supersymmetry is the odd-odd-odd bracket. Furthermore, due to polarization it is enough to ask for the vanishing of 
\begin{equation}
\begin{aligned}
\mathcal{L}_{\{\vec\xi,\vec\xi\}}\vec\xi
=\ & \{\vec\xi,\vec\xi\}^+\mathcal{L}_{K_{(+)}}\vec\xi 
	+\{\vec\xi,\vec\xi\}^-\mathcal{L}_{K_{(-)}}\vec\xi 
	+\sum_i\{\vec\xi,\vec\xi\}^i \mathcal{L}_{K_{(i)}}\vec\xi \\
&	+\sum_i\{\vec\xi,\vec\xi\}^{i^*} \mathcal{L}_{K_{(i^*)}}  \vec\xi
	+\frac{1}{2}\sum_{ij}\{\vec\xi,\vec\xi\}^{ij} \mathcal{L}_{K_{(ij)}}\vec\xi
\end{aligned}
\end{equation}
for all $\vec\xi\in\mathcal{K}_1$.

In our situation we use the notations from Section \ref{sec:geoalg}, in particular ({\ref{50}), (\ref{proj:W}), (\ref{proj:V*}), and (\ref{proj:so}), such that the vanishing of $\mathcal{L}_{\{\vec\xi,\vec\xi\}}\vec\xi$ for $\vec\xi=\vec\xi(\xi_1^0,\xi_2^0)$ is 
\begin{align*}
0  
%=\ &
%	-C\big(\xi_2^0,\Gamma_+\xi_2^0\big)\, \vec\xi\big(\bar c\xi^0_1, d\xi_2^0\big) \\
%&	-\frac{1}{8}\sum_{ij} C\big(\xi^0_2, \Gamma_+ \big(s_{d,\bar c}(e_j)\Gamma_i
%	+\Gamma_{i}s_{\bar c,d}(e_j)\big)\xi^0_2\big)
%	\,\vec\xi\big(\Gamma_{ij}\xi_1^0,\Gamma_{ij}\xi_2^0\big)\\
%&	+\sum_i C\big(\xi_1^0,\Gamma_i\xi_2^0\big)\, \vec\xi\big(\Gamma_+s_{\bar c,d}(e_i)\xi_2^0,0\big) \\
%&	-\sum_i C\big(\xi_1^0,s_{\bar c,d}(B^{-1}(e_i))\xi_2^0\big)
%	\,\vec\xi\big(\Gamma_+B(e_i)\xi_2^0,0\big)\\
=\ &
	-C\big(\xi_2^0,\Gamma_+\xi_2^0\big)
		\, \vec\xi\big(\bar c\xi^0_1, d\xi_2^0\big) 
	+\frac{i\lambda_6}{4}\sum_{ij=6}^9 C\big(\xi^0_2, \Gamma_+ \Gamma_{125ij}\xi^0_2\big)
		\,\vec\xi\big(\Gamma_{ij}\xi_1^0,\Gamma_{ij}\xi_2^0\big)\\
&	+\frac{i\lambda_1}{2}C\big(\xi^0_2, \Gamma_+ \Gamma_5 \xi^0_2\big)
		\,\vec\xi\big(\Gamma_{12}\xi_1^0,\Gamma_{12}\xi_2^0\big)
	+\frac{i\lambda_3}{2} C\big(\xi^0_2, \Gamma_+ \Gamma_{12345} \xi^0_2\big)
		\,\vec\xi\big(\Gamma_{34}\xi_1^0,\Gamma_{34}\xi_2^0\big)\\
&	+i\sum_i \lambda_i C\big(\xi_1^0,\Gamma_i\xi_2^0\big)
		\, \vec\xi\big(\Gamma_+\Gamma_{125}\Gamma_i\xi_2^0,0\big) 
	-i\sum_i \lambda_i C\big(\xi_1^0,\Gamma_{125}\Gamma_i\xi_2^0\big)
		\,\vec\xi\big(\Gamma_+ \Gamma_i\xi_2^0,0\big)\,.
\end{align*}
Along with $\vec\xi(\xi_1^0,\xi_2^0)=0\Leftrightarrow\xi_1^0=\xi_2^0=0$ this yields the following two equations for the two constant spinors: 
\begin{align*}
0
=\ & -\alpha_+'C\big(\xi_2^0,\Gamma_+\xi_2^0\big)\, \Gamma_{125}\xi_2^0  
    +\frac{\alpha_++\alpha'_+}{4}\sum_{ij=6}^9 
    	C\big(\xi^0_2, \Gamma_+\Gamma_{125ij}\xi^0_2\big)\Gamma_{ij}\xi_2^0\\
&	+\frac{\alpha_--\alpha'_+}{2}C\big(\xi^0_2, \Gamma_+ \Gamma_5 \xi^0_2\big)\Gamma_{12}\xi_2^0
	+\frac{\alpha_-+\alpha'_+}{2}C\big(\xi^0_2, \Gamma_+ \Gamma_{12345} \xi^0_2\big)\Gamma_{34}\xi_2^0
\end{align*}
and
\begin{align*}
0
%=\ &
%	-C\big(\xi_2^0,\Gamma_+\xi_2^0\big)(\alpha_+X_{1234}^++\alpha_-X^-_{1234}\Gamma_{125}\xi^0_1 
%	+\frac{i\lambda_1}{2}C\big(\xi^0_2, \Gamma_+ \Gamma_5 \xi^0_2\big)\Gamma_{12}\xi_1^0\\
%&	+\frac{i\lambda_3}{2} C\big(\xi^0_2, \Gamma_+ \Gamma_{12345} \xi^0_2\big)\Gamma_{34}\xi_1^0
%	+\frac{i\lambda_6}{4}\sum_{ij=6}^9 C\big(\xi^0_2, \Gamma_+ \Gamma_{125ij}\xi^0_2\big)\Gamma_{ij}%\xi_1^0\\
%&	+i\sum_i \lambda_i C\big(\xi_1^0,\Gamma_i\xi_2^0\big)\, \Gamma_{125}\Gamma_i\xi_2^0 
%	-i\sum_i \lambda_i C\big(\xi_1^0,\Gamma_{125}\Gamma_i\xi_2^0\big)\Gamma_i\xi_2^0 \\
=\ &
	-\frac{\alpha_+}{2} C\big(\xi_2^0,\Gamma_+\xi_2^0\big)(\Gamma_{125}-\Gamma_{345})\xi^0_1  
	-\frac{\alpha_-}{2} C\big(\xi_2^0,\Gamma_+\xi_2^0\big)(\Gamma_{125}+\Gamma_{345})\xi^0_1 \\
&	+\frac{\alpha_--\alpha_+'}{2}C\big(\xi^0_2, \Gamma_+ \Gamma_5 \xi^0_2\big)\Gamma_{12}\xi_1^0
	+\frac{\alpha_-+\alpha_+'}{2} C\big(\xi^0_2, \Gamma_+ \Gamma_{12345} \xi^0_2\big)\Gamma_{34}\xi_1^0\\
&	+\frac{\alpha_++\alpha_+'}{4}\sum_{ij=6}^9 C\big(\xi^0_2, \Gamma_+ \Gamma_{125ij}\xi^0_2\big)	
		\Gamma_{ij}\xi_1^0\\
&	+(\alpha_--\alpha_+')\sum_{i=1}^2  \big(
		C\big(\xi_1^0,\Gamma_i\xi_2^0\big)\, \Gamma_{125}\Gamma_i\xi_2^0
		- C\big(\xi_1^0,\Gamma_{125}\Gamma_i\xi_2^0\big)\Gamma_i\xi_2^0 \big)\\
&	+(\alpha_-+\alpha_+')\sum_{i=3}^4 \big( 
		C\big(\xi_1^0,\Gamma_i\xi_2^0\big)\, \Gamma_{125}\Gamma_i\xi_2^0 
		- C\big(\xi_1^0,\Gamma_{125}\Gamma_i\xi_2^0\big)\Gamma_i\xi_2^0	\big)\\
&	-(\alpha_+-\alpha_+')\big( 
		C\big(\xi_1^0,\Gamma_5\xi_2^0\big)\, \Gamma_{12}\xi_2^0 
		- C\big(\xi_1^0,\Gamma_{12}\xi_2^0\big)\Gamma_5\xi_2^0\big)\\
&	+(\alpha_++\alpha_+')\sum_{i=6}^9 \big(
		C\big(\xi_1^0,\Gamma_i\xi_2^0\big)\, \Gamma_{125}\Gamma_i\xi_2^0 
		- C\big(\xi_1^0,\Gamma_{125}\Gamma_i\xi_2^0\big)\Gamma_i\xi_2^0\big) \,.
\end{align*}
We may rewrite this in terms of $V=\RR^9$ by considering $\xi_1^0,\xi_2^0\in S(V)$ and using the conventions from Section \ref{sec:conv}
\begin{equation}\label{ooo2}
\begin{aligned}
0
=\ & -\alpha_+'C_V\big(\xi_2^0,\xi_2^0\big)\, \gamma_{125}\xi_2^0  
    +\frac{\alpha_++\alpha'_+}{4}\sum_{ij=6}^9 
    	C_V\big(\xi^0_2,\gamma_{125ij}\xi^0_2\big)\gamma_{ij}\xi_2^0\\
&	+\frac{\alpha_--\alpha'_+}{2}C_V\big(\xi^0_2, \gamma_5 \xi^0_2\big)\gamma_{12}\xi_2^0
	+\frac{\alpha_-+\alpha'_+}{2}C_V\big(\xi^0_2,  \gamma_{12345} \xi^0_2\big)\gamma_{34}\xi_2^0
\end{aligned}
\end{equation}
and
\begin{align}
0
=\ &
	\frac{\alpha_++\alpha_-}{2} C_V\big(\xi_2^0,\xi_2^0\big)\gamma_{125}\xi^0_1  
	-\frac{\alpha_+-\alpha_-}{2} C_V\big(\xi_2^0,\xi_2^0\big)\gamma_{345}\xi^0_1 \nonumber\\
&	+\frac{\alpha_--\alpha_+'}{2}C_V\big(\xi^0_2, \gamma_5 \xi^0_2\big)\gamma_{12}\xi_1^0
	+\frac{\alpha_-+\alpha_+'}{2} C_V\big(\xi^0_2, \gamma_{12345} \xi^0_2\big)\gamma_{34}\xi_1^0\nonumber\\
&	+\frac{\alpha_++\alpha_+'}{4}\sum_{ij=6}^9 C_V\big(\xi^0_2, \gamma_{125ij}\xi^0_2\big)	
		\gamma_{ij}\xi_1^0\nonumber\\
&	-2(\alpha_--\alpha_+')\sum_{i=1}^2  \big(
		C_V\big(\xi_1^0,\gamma_i\xi_2^0\big)\, \delta_{i[1}\gamma_{2]5}\xi_2^0
		- C_V\big(\xi_1^0,\delta_{i[1}\gamma_{2]5}\xi_2^0\big)\gamma_i\xi_2^0 \big)\label{ooo1}\\
&	+(\alpha_-+\alpha_+')\sum_{i=3}^4 \big( 
		C_V\big(\xi_1^0,\gamma_i\xi_2^0\big)\, \gamma_{125i}\xi_2^0 
		- C_V\big(\xi_1^0,\gamma_{125i}\xi_2^0\big)\gamma_i\xi_2^0	\big)\nonumber\\
&	-(\alpha_+-\alpha_+')\big( 
		C_V\big(\xi_1^0,\gamma_5\xi_2^0\big)\, \gamma_{12}\xi_2^0 
		- C_V\big(\xi_1^0,\gamma_{12}\xi_2^0\big)\gamma_5\xi_2^0\big)\nonumber\\
&	+(\alpha_++\alpha_+')\sum_{i=6}^9 \big(
		C_V\big(\xi_1^0,\gamma_i\xi_2^0\big)\, \gamma_{125i}\xi_2^0 
		- C_V\big(\xi_1^0,\gamma_{125i}\xi_2^0\big)\gamma_i\xi_2^0\big) \,.
\end{align}

A direct computation shows that (\ref{ooo1}) and (\ref{ooo2}) are obtained only for $\alpha_+=-3\alpha_+'$.

\begin{thm}\label{thm:susy}
The geometric superalgebra of the indecomposable Cahen-Wallach space $M_B$ according to Theorem \ref{thm:geoalg} yields non restricted geometric supersymmetry if and only if 
\[
B=-{\rm diag}\left(
\left(\alpha_--\alpha_+'\right)^2\mathbbm{1}_2,\,
\left(\alpha_-+\alpha_+'\right)^2\mathbbm{1}_2,\,
16\alpha_+'^2\mathbbm{1}_1,\,
4\alpha_+'^2\mathbbm{1}_4  
\right)
\]
and $(\bar c,d)$ given by
%\,\footnote{Strictly speaking, we may add an arbitrary summand $\beta'X^-_{1234}\Gamma_{125}$ to $d$ without changing the result because $X^-_{1234}\xi_2^0=0$.}
\[
\bar c= \big(-3\alpha_+' X^+_{1234}+\alpha_- X^-_{1234}\big)\Gamma_{125},\quad 
d = \big(\alpha_+' X^+_{1234}+\alpha'_- X^-_{1234}\big)\Gamma_{125}  \,.
\]
\end{thm}

\section{The moduli space of geometric supersymmetry}\label{sec:7}

\subsection{The moduli space of geometric superalgebras and supersymmetries}

The moduli space of geometric superalgebras according to Theorem \ref{thm:geoalg} is naturally parameterized by 
\[
(\alpha_-,\alpha_+',\alpha_+,\alpha'_-)
	\in \RR^4\setminus \left\{
	(\pm\alpha_+',\alpha_+',\alpha_+,\alpha'_-),(\alpha_-,\alpha_+',\pm\alpha_+',\alpha'_-)\,|\,
			\alpha_\pm,\alpha_\pm'\in\RR
	\right\}.
\]
If we exclude the euclidean configuration and divide out the isometries defined by the action of positive(!)\ scalars we end up with a 3-sphere with several points removed. 

If we in addition divide out the isometries defined by 
\[(\alpha_-,\alpha_+',\alpha_+,\alpha'_-)\sim(-\alpha_-,\alpha_+',\alpha_+,\alpha_-')\]
we are left with a closed half sphere that can be parameterized by $\alpha_+,\alpha_+',\alpha_-'$ via $\alpha_-=\sqrt{1-\alpha_+^2-{\alpha_+'}^2-{\alpha_-'}^2}$. 
At the end we obtain a subset of the closed 3-ball $D^3$
\begin{equation}\label{mod-1}
\mathcal{C}
	 =  D^3\setminus\left\{\left(\alpha_+',\alpha_+,\alpha_-'\right)\,\middle|\, \alpha_+=\pm\alpha_+' 
	 		\text{ or } \alpha_+^2+2\alpha'_+{}^2+\alpha_-'{}^2=1\right\}\,,
\end{equation}
see Figure \ref{fig:moduli}. 

\begin{figure}\caption{The moduli space $\widehat{\mathcal{C}}$ of geometric superalgebras projected to the $\alpha_+\alpha_+'$-plane.}
	\label{fig:moduli}
	\includegraphics[scale=0.6]{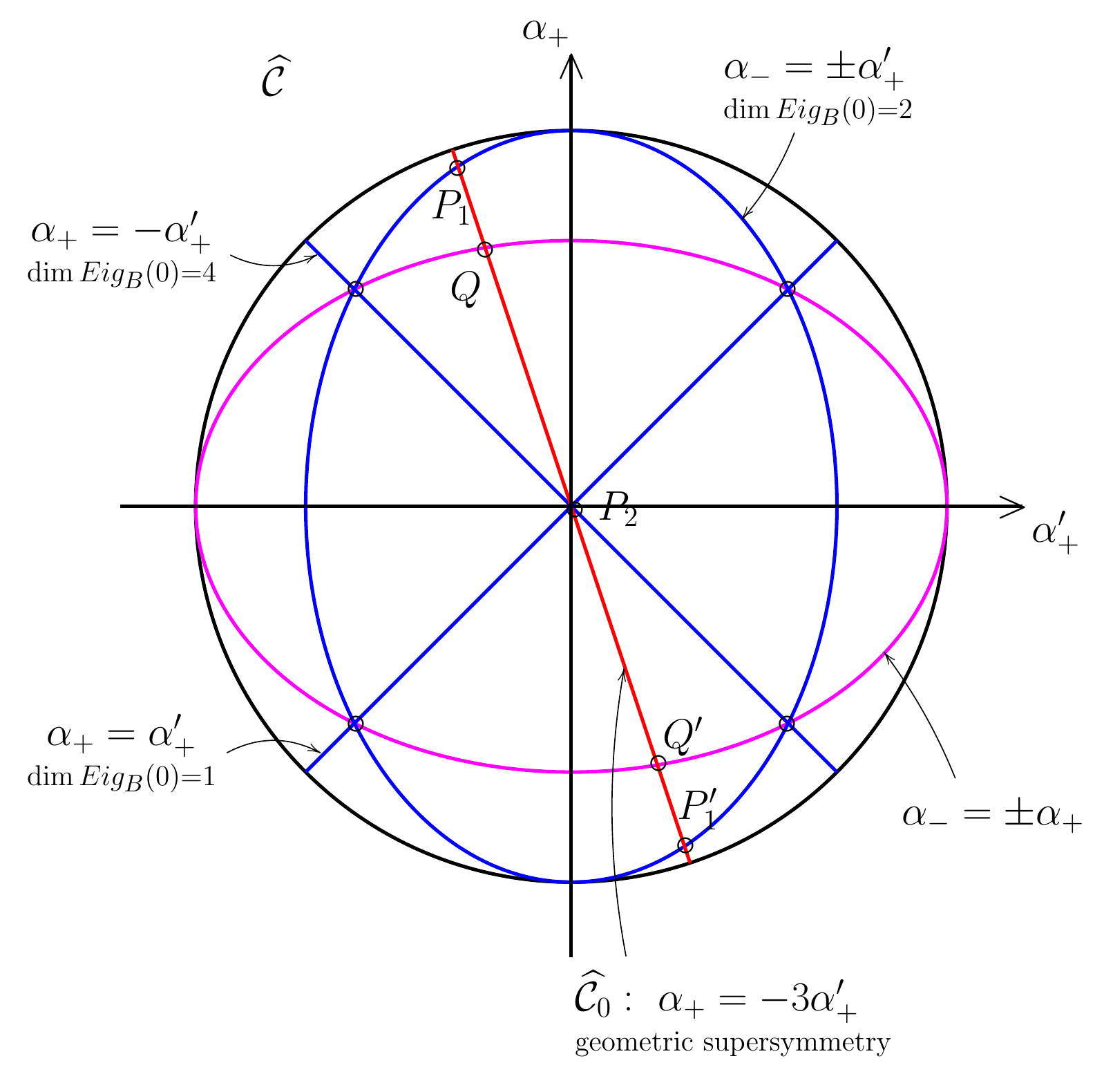}
\end{figure}

The excluded parameters in (\ref{mod-1}) yield decomposable spaces: On the 2-discs (the diagonals in Figure \ref{fig:moduli}) the map $B$ admits either one zero eigenvalue ($\alpha_+=\alpha_+'$) or four zero eigenvalues ($\alpha_+=-\alpha_+'$) and on the ellipsoid that is defined by $\alpha_-^2={\alpha_+'}^2\Leftrightarrow \alpha_+^2+2{\alpha_+'}^2+{\alpha_-'}^2=1$ (indicated by the ellipse through $P_1$ and $P_1'$ in Figure \ref{fig:moduli}) it admits two zero eigenvalues.
In the intersection points of the sets the multiplicities have to be summed up.
\begin{prop}
The compactified moduli space of geometric superalgebras according to Theorem \ref{thm:geoalg} is $\widehat{\mathcal{C}}=D^3$, see (\ref{mod-1}). It is obtained by adding the decomposable, non-euclidean configurations indicated by the diagonal 2-discs and the ellipsoid through $P_1$ in Figure \ref{fig:moduli}. 
\end{prop}

\begin{rem}
If we divide out the remaining isometries defined by the antipodal map the resulting space is an $\RR P^2$-cone over the base point $P_2$ and the decomposable spaces are associated to special sections: Two $S^1$ cones over $P_2$ as well as one ''ellipsoid''. We will not use this description below. In Figure \ref{fig:moduli} such isometric configurations will be denoted by non-primed/primed pairs.
\end{rem}

\begin{prop}
The compactified moduli space $\widehat{\mathcal{C}_0}$ of geometric supersymmetries according to Theorem \ref{thm:susy} is parameterized by
the 2-disc defined by $\alpha_+=-3\alpha_+'$ in $B^3$. It is is indicated by the red line in Figure \ref{fig:moduli}. 
\end{prop}

\subsection{The singluar points as $N$-extended supersymmetries}

\begin{rem}\label{rem:N-ext}
	Up to now, we considered algebras that have been constructed from spinor bundles associated to irreducible Clifford modules $S_0$. We may instead consider a spinor bundle that is associated to a reducible module of the form $S_0\otimes\CC^N$.  Then the results from Propositions \ref{prop:maps} and \ref{prop:flat} remain formally the same if we consider the parameters to take their values in $\cl(W)\otimes\mathfrak{gl}_N\CC$. We provide two examples in (\ref{pair-i}) and (\ref{pair-ii}) below.
	
	If a geometric superalgebra in the sense of Section \ref{sec:geoalg1} is based on a reducible spinor bundle of the above type then we call it $N$-extended. 
\end{rem}
There are two configurations of parameters in the compactified moduli $\widehat{\mathcal{C}}_0$ of geometric supersymmetries that yield decomposable spaces. These are 
\begin{enumerate}
\item[{\bf i.}] 
the ellipse through the points $P_1$ and $P'_1$ defined by $\alpha_-=\pm\alpha_+$ and 
\item[{\bf ii.}]
the line through $P_2$ defined by $\alpha_+'=0$ 
\end{enumerate}
that are contained in the disc $\widehat{\mathcal{C}}_0$ in Figure \ref{fig:moduli}.
If we specify the parameter $\alpha_-'$ to  $\alpha_-'=0$ we exactly reach the mentioned points:  
\begin{enumerate}
\item[$P_2$] with $B=-\alpha_-^2{\rm diag}\left(\mathbbm{1}_4,\mathds{O}_5 \right)\,.$
 
In this case we have $c=\alpha_- X^-_{1234}\Gamma_{125}$ and $d=0$.
Therefore, the action on $\sigma_-X_{1234}^+S$ is trivial, and we may further reduce $\mathcal{K}_1$ to 
\[
X_{1234}^-S_{11}^-\oplus X_{1234}^+S_{11}^+\subset S_{11}^-\oplus S_{11}^+= S_{11}\,.
\]
\item [$P_1$] with $B=-4\alpha_-^2{\rm diag}\left(4\mathbbm{1}_1,\mathbbm{1}_6,\mathds{O}_2\right).$

In this case a further reduction is not possible.
\end{enumerate}
We emphasize the fact, that in both cases {\bf i.} and {\bf ii.} conditions (\ref{ooo2}) and (\ref{ooo1}) are independent of the Killing vector fields associated to the zero eigenvalues of $B$.

The geometric supersymmetries on the decomposable eleven dimensional spaces that are associated to the singular points of the moduli space can be interpreted as $N$-extended geometric supersymmetries in lower dimensions $D<11$, at least if we restrict the even part $\mathcal{K}_0$ in a suitable way, i.e.\ to $i,i^*\in\{1,\ldots,D-2\}$:
 
The singular point $P_2$ can be associated to restricted $\nu=\sfrac{1}{2}$, 4-extended geometric supersymmetry in six dimensions. The ingredients are as follows

\begin{itemize}[leftmargin=4.5ex]
\item
	The $D=6$ Cahen-Wallach space $M_6$ associated to $B=-\beta^2\mathbbm{1}_4$.
\item
	The spinor bundle $S=S(M_6)\otimes \CC^4$, 
\item
	The bilinear form $C=C_6\otimes C_5$ with $C_5$ being the charge conjugation on $S_5=\CC^4$. 
\item
	The non-flat connection according to Proposition \ref{prop:maps} and Remark \ref{rem:N-ext} that is defined by  
\begin{equation}\label{pair-i}
\bar c=\beta X^-_{1234}\Gamma^{(6)}_{12} \otimes T,\quad 
d=0
\end{equation}
with $T$ being some vector in $\cl_1(\RR^5)$ with $T^2=-\mathbbm{1}$. 
\item
	The even part is defined by the Killing vector fields of $M_6$.
\item
	The odd part $\mathcal{K}_1$ is defined by 
\[
\left(X^-_{1234}S_6^-\oplus X^+_{1234}S_6^+\right)\otimes\CC^4= \Pi^+S_6 \otimes\CC^4\,.
\]
The space we just described is exactly the $D=6,N=4$ supergravity background discussed in \cite{Meessen2002}.
\end{itemize}

The singular point $P_1$ can be associated to non-restricted, i.e.\ $\nu=\sfrac{3}{4}$, 2-extended geometric supersymmetry in nine dimensions. Here the correspondence is as follows.

\begin{itemize}[leftmargin=4.5ex]
\item 
	The $D=9$ Cahen-Wallach space $M_9$ associated to $B=-4\alpha^2{\rm diag}\left(4 \mathbbm{1}_1,\mathbbm{1}_6\right)$. 
\item
	The spinor bundle $S=S(M_9)\otimes \CC^2$.
\item
	The bilinear form $C=C_9\otimes \sigma_1$.
\item
	The non-flat connection according to Proposition \ref{prop:maps} and Remark \ref{rem:N-ext} that is defined by 
\begin{equation}\label{pair-ii}
c=-\alpha\Gamma^{(9)}_{1}\otimes\mathbbm{1}_2 + 2\alpha\Gamma^{(9)}_{123}\otimes i\sigma_3,\quad 
d= \frac{\alpha}{2}\Gamma^{(9)}_1\otimes\mathbbm{1}_2 +\frac{\alpha}{2}\Gamma^{(9)}_{123}\otimes i\sigma_3\,.
\end{equation}
\item
	The even part is defined by the Killing vector fields of $M_9$.
\item
	The odd part $\mathcal{K}_1$ is then defined by 
\[
\big(S_7\oplus X_{23}^-S_7) \oplus \big(S_7\oplus X_{23}^+S_7\big) 
	\subset (S_2\otimes S_7)\oplus(S_2\otimes S_7) 
	= S_9\otimes \CC^2\,.
\]
\end{itemize}

\begin{rem}\phantomsection\label{rem:D6N4}
	We will take a closer look at Examples \ref{ex:maximal} and \ref{ex:maximal-cont}. 
\begin{enumerate}[leftmargin=4.5ex]
\item
	In the moduli space $\widehat{\mathcal{C}}$ we may consider the ellipsoid defined by $\alpha_+=\pm\alpha_-$ (this is indicated by the ellipse in Figure \ref{fig:moduli} running through ${Q}$ and ${Q'}$). Here the symmetric map $B$ has exactly two different eigenvalues and is given by 
	\[
	B=-\begin{pmatrix}(\alpha_+-\alpha_+')^2\mathbbm{1}_3&\\&(\alpha_++\alpha_+')^2\mathbbm{1}_{6}\end{pmatrix}\,.
	\]
	A special situation occurs, when we furthermore fix the (almost free) parameter $\alpha_-'$ to $\alpha_-'=\pm\alpha_+'$, where we choose the same sign as above.
	The result is an ellipse that is the intersection of the above ellipsoid with the disc obtained by $\alpha_-'=\pm\alpha_+'$.   
	In this situation the connection is determined by\footnote{This holds for the upper sign, the lower sign interchanges the role of $(12)$ and $(34)$.}
	$(\bar c,d)=\left(\alpha_+\Gamma_{125}^{(11)}, \alpha_+'\Gamma_{125}^{(11)}\right)$ and the related geometric superalgebra is flat. Therefore, the restriction of the spinors is no longer necessary, see Proposition \ref{prop:flat} and \eqref{q}.
  	
	In particular, there is a pair of points $P_0,P'_0$ defined by $\alpha_+=-3\alpha_+'$ for which non-restricted $\nu=\sfrac{1}{1}$ geometric supersymmetry is achieved. 
	In Figure \ref{fig:moduli} they are given by the intersection of the mentioned ellipse with the disc $\widehat{\mathcal{C}}_0$. One of the two points lie above and the other below the $\alpha_+\alpha_+'$-plane and the full configuration in $S^3$ is given by $(\alpha_-,\alpha_+',\alpha_+,\alpha_-')=\frac{1}{2\sqrt{5}}(\pm3,\mp1,\pm3,\mp1)$. 
	The connection in this case yields the unique supergravity background found by \cite{CheKo84,FigPapado1}. \\
\item
	If we restrict the odd part of the $\nu=\sfrac{1}{1}$ geometric superalgebra 
%on the ellipsoid $\alpha_=\pm\alpha_+$ 
	to $X_{1234}^-S_{11}^-\oplus X^+_{1234}S_{11}$ and consider $\alpha_+'=0$ we see the following feature: 	
	Although the analog to (\ref{ooo2}) and (\ref{ooo1}) is not obtained for the full summation $1,\ldots,9$, it is obtained for a summation $1,\ldots,4$. 	
	Therefore, if we again restrict the even part in a suitable way, i.e.\ to $1,\ldots,4$, we get a super Lie algebra that can be interpreted as the same restricted $\nu=\sfrac{1}{2}$, 4-extended geometric supersymmetry as before, but with $(\bar c,d)=\big(\beta\Gamma^{(6)}_{12}\otimes T,0\big)$ instead. 
	The main differences to the interpretation in $P_2$ is that in this case the eleven dimensional oxidation is flat, indecomposable, and defines a geometric superalgebra only, instead of geometric supersymmetry.
\end{enumerate}
\end{rem}

Although we have been very brief in the description of the two singular points, we hope that the reader is well prepared to handle these example by using the preliminaries provided in this text.

\appendix

\section{Clifford conventions}\label{sec:conv}

Our convention for the Clifford algebra $\cl(V,g)$ of a pseudo-Riemannian space $(V,g)$ is given by
\begin{equation}
	vw+wv=-2g(v,w)
\end{equation}
for all $v,w\in V$. This yields the spin representation $\Gamma: \mathfrak{so}(V,g)\to\cl(V,g)$ as follows. We consider a basis $\{e_i\}$ of $V$ and denote its image in $\cl(V,g)$ by $\{\Gamma_i\}$. We denote the basis of $\mathfrak{so}(V,g)$ by $E_{ij}$ with $E_{ij}(e_k)=2g_{k[j}e_{i]}$ such that for $A=\frac{1}{2}\sum_{ij}A_{ij}E^{ij}$ and $v=\sum_kv^ke_k$ we have $A(v)=\sum_{ij}A^i{}_{j}v^je_i$. Then the image of $E_{ij}$ under the spin representation $\Gamma$ is given by 
\begin{equation}\label{spinrepr}
\Gamma(E_{ij})=-\frac{1}{2}\Gamma_{ij} \quad \text{or}\quad \Gamma(A)=-\frac{1}{4}A_{ij}\Gamma^{ij}\,.
\end{equation}
In particular, the metric isomorphism $\Lambda^2V\simeq \mathfrak{so}(V,g)$ given by $e_i\wedge e_j \leftrightarrow E_{ij}$ along with (\ref{spinrepr}) yields (\ref{vstar}). 

For any vector field $Y\in\mathfrak{X}(M)$ on a pseudo Riemannian spin manifold $(M,g)$ the Levi-Civita connection $\nabla$ on $M$ yields a linear map $\nabla Y:\mathfrak{X}(M)\to\mathfrak{X}(M)$ by $\nabla Y(X):=\nabla_X Y$. 
Using a local basis $\{\partial_i\}$ of $M$ we have $(\nabla Y)_{ij}=\nabla_j Y_i$.

If we denote the Christoffel symbols by $\Gamma_{ij}^k$ such that 
$\nabla_iY=\partial_iY+\sum_{j k}\Gamma_{ij}^k Y^j\partial_k$  the connection form $\omega$ with $\nabla_iY =\partial_iY +\omega_iY=\partial_iY+\frac{1}{2}\sum_{jk}(\omega_i)_{jk}E^{jk}(Y)$ is given by $(\omega_i)_{jk}=\Gamma_{ik;j}$ Therefore, the induced connection on the spinor bundle is given by
\begin{equation}
\nabla_i\psi
 =\partial_i\psi-\frac{1}{4}\sum_{jk}(\omega_{i})_{jk}\Gamma^{jk}\psi
 =\partial_i\psi+\frac{1}{4}\sum_{jk}\Gamma_{ij;k}\Gamma^{jk}\,.
\end{equation} 
Moreover, the spinorial Lie derivative that is defined by (\ref{spin-Lie}) reads as 
\begin{equation}
\mathcal{L}_K\psi=\nabla_K\psi-\Gamma(\nabla K)\psi
= \nabla_k\psi+\frac{1}{4}\sum_{ij}(\nabla K)_{ij}\Gamma^{ij}
= \nabla_k\psi-\frac{1}{4}\sum_{ij}\nabla_{i}K_j\Gamma^{ij}\,.
\end{equation}

We use the following explicit Clifford representations in dimension nine and eleven for the calculations in Section \ref{sec:susy}.

Consider matrices $L_a$ for $1\leq a\leq 7$ that are defined as matrix representation of left multiplication by imaginary octonions:
{\renewcommand{\arraystretch}{0}\renewcommand{\tabcolsep}{0.1pt}\small
\begin{align*}
L_1&=\left(\text{\begin{tabular}{cccccccc}
 &--1&  &  &  &  &  & \\
{ }1&  &  &  &  &  &  & \\
 &  &  &--1&  &  &  & \\
 &  &{ }1 &  &  &  &  & \\
 &  &  &  &  &--1&  & \\
 &  &  &  &{ }1 &  &  & \\
 &  &  &  &  &  &  &{ }1 \\
 &  &  &  &  &  &--1&
 \end{tabular}}\right)
,&
L_2&=\left(\text{\begin{tabular}{cccccccc}
 &  &--1&  &  &  &  & \\
 &  &  &{ }1 &  &  &  & \\
{ }1&  &  &  &  &  &  & \\
 &--1&  &  &  &  &  & \\
 &  &  &  &  &  &--1& \\
 &  &  &  &  &  &  &--1 \\
 &  &  &  &{ }1 &  &  & \\
 &  &  &  &  &{ }1 &  &
 \end{tabular}}\right)
,&
L_3&=\left(\text{\begin{tabular}{cccccccc}
 &  &  &--1&  &  &  & \\
 &  &--1&  &  &  &  & \\
 &{ }1 &  &  &  &  &  & \\
{ }1&  &  &  &  &  &  & \\
 &  &  &  &  &  &  &--1\\
 &  &  &  &  &  &{ }1 &  \\
 &  &  &  &  &--1&  & \\
 &  &  &  &{ }1 &  &  &
 \end{tabular}}\right)
,\\
L_4&=\left(\text{\begin{tabular}{cccccccc}
 &  &  &  &--1&  &  & \\
 &  &  &  &  &{ }1 &  & \\
 &  &  &  &  &  &{ }1 & \\
 &  &  &  &  &  &  &{ }1 \\
{ }1&  &  &  &  &  &  & \\
 &--1&  &  &  &  &  &  \\
 &  &--1&  &  &  &  & \\
 &  &  &--1&  &  &  &
 \end{tabular}}\right)
,&
L_5&=\left(\text{\begin{tabular}{cccccccc}
 &  &  &  &  &--1&  & \\
 &  &  &  &--1&  &  & \\
 &  &  &  &  &  &  &{ }1\\
 &  &  &  &  &  &--1& \\
 &{ }1 &  &  &  &  &  & \\
{ }1&  &  &  &  &  &  & \\
 &  &  &{ }1 &  &  &  & \\
 &  &--1&  &  &  &  &
 \end{tabular}}\right)
,&
L_6&=\left(\text{\begin{tabular}{cccccccc}
 &  &  &  &  &  &--1& \\
 &  &  &  &  &  &  &--1\\
 &  &  &  &--1&  &  & \\
 &  &  &  &  &{ }1 &  & \\
 &  &1 &  &  &  &  & \\
 &  &  &--1&  &  &  & \\
{ }1&  &  &  &  &  &  & \\
 &{ }1 &  &  &  &  &  &
 \end{tabular}}\right)
,\\
L_7&=\left(\text{\begin{tabular}{cccccccc}
 &  &  &  &  &  &  &--1\\
 &  &  &  &  &  &{ }1 & \\
 &  &  &  &  &--1&  & \\
 &  &  &  &--1&  &  & \\
 &  &  &{ }1 &  &  &  & \\
 &  &{ }1 &  &  &  &  & \\
 &--1&  &  &  &  &  & \\
{ }1&  &  &  &  &  &  &
 \end{tabular}}\right)\,. 
\end{align*}}
Starting from this we define $\gamma$-matrices $\{\gamma_i\}_{i\in\{1,\ldots,9\}}$ for the euclidean space $V=\RR^9$ by 
\begin{align*}
\gamma_a:=\sigma_1\otimes L_a\,,\quad
\gamma_8=-i\sigma_2\otimes \mathbbm{1}\,,\quad
\gamma_9:=-i\sigma_3\otimes\mathbbm{1}\,.
\end{align*}
The charge conjugation matrix for $V$ obeys $ \gamma_i^tC_V =C_V\gamma_i $ and is given by
\[
C_V=\sigma_3\otimes \mathbbm{1}\,.
\]
From this we get the $\gamma$-matrices $\{\Gamma_\mu\}_{\mu\in\{+,-,i\}}$ on $W=\RR^{1,10}$ by the procedure described in Section \ref{sec:prel}. In particular, the charge conjugation obeys $\Gamma_\mu^tC_W=C_W\Gamma_\mu$ and is given by
\[
C_W=\sigma_2\otimes C_V \,.
\]
Here as well as in the main text the matrices 
$\sigma_1=\begin{pmatrix}0&1\\1&0\end{pmatrix}$, 
$\sigma_2=\begin{pmatrix}0&-i\\i&0\end{pmatrix}$, and 
$\sigma_3=\begin{pmatrix}1&0\\0&-1\end{pmatrix}$, 
that obey $\sigma_i\sigma_j=i\sum_{k=1}^3\epsilon_{ijk}\sigma_k$ 
denote the Pauli-matrices.

%\bibliographystyle{alpha}\addcontentsline{toc}{section}{References}
%\bibliography{bibliothek}

\end{document}